\documentclass[11pt,a4paper]{amsart}
\usepackage{amscd,amsfonts,amsmath,amsthm,amssymb,verbatim,mathrsfs}
\usepackage[dvips]{graphicx}
\usepackage{tikz}
\usetikzlibrary{shapes,backgrounds}
\usetikzlibrary{calc}
\usetikzlibrary{matrix}
\usepackage{graphics,hyperref}
\usepackage{pstricks}
\usepackage{enumitem}
\usepackage{psfrag}
\usepackage{setspace}
\usepackage{latexsym}
\usepackage{caption}
\usepackage{pstcol,pst-plot,pst-3d}
\usepackage[normalem]{ulem} 

\newpsstyle{fatline}{linewidth=1.5pt}
\newpsstyle{fyp}{fillstyle=solid,fillcolor=verylight}
\definecolor{verylight}{gray}{0.97}
\definecolor{light}{gray}{0.9}
\definecolor{medium}{gray}{0.85}

\def\frk{\frak}               

\def\Phi{{\frk n}}
\def\Phi{{\frk N}}
\def\opn#1#2{\def#1{\operatorname{#2}}} 
\opn\chara{char} \opn\length{\ell} \opn\pd{pd} \opn\rk{rk}
\opn\projdim{proj\,dim} \opn\injdim{inj\,dim} \opn\rank{rank}
\opn\depth{depth} \opn\grade{grade} \opn\height{height}
\opn\embdim{emb\,dim} \opn\codim{codim} \opn\sgn{sgn}

\opn\Tr{Tr} \opn\bigrank{big\,rank}
\opn\superheight{superheight}\opn\lcm{lcm}
\opn\trdeg{tr\,deg}
\opn\reg{reg} \opn\lreg{lreg} \opn\ini{in} \opn\lpd{lpd}
\opn\size{size}\opn\bigsize{bigsize}
\opn\cosize{cosize}\opn\bigcosize{bigcosize}
\opn\sdepth{sdepth}\opn\sreg{sreg}
\opn\link{link}\opn\fdepth{fdepth}
\usepackage{tikz}
\usetikzlibrary{shapes,backgrounds}
\usetikzlibrary{calc}
\usetikzlibrary{matrix}
\usetikzlibrary{arrows}
\usepackage{graphics,hyperref}
\usepackage{enumitem}
\usepackage{psfrag}
\opn\div{div} \opn\Div{Div} \opn\cl{cl} \opn\Cl{Cl} \opn\Cor{Cor}
\opn\Spec{Spec} \opn\Supp{Supp} \opn\supp{supp} \opn\Sing{Sing}
\opn\Ass{Ass} \opn\Min{Min}\opn\Mon{Mon} \opn\dstab{dstab} \opn\astab{astab}
\opn\Ann{Ann} \opn\Rad{Rad} \opn\Soc{Soc} \opn\Gr{Gr}
\opn\Im{Im} \opn\Ker{Ker} \opn\Coker{Coker} \opn\Am{Am}
\opn\Hom{Hom} \opn\Tor{Tor} \opn\Ext{Ext} \opn\End{End}
\opn\Aut{Aut} \opn\id{id} \opn\span{span}

\opn\nat{nat}
\opn\pff{pf}
\opn\Pf{Pf} \opn\GL{GL} \opn\SL{SL} \opn\mod{mod} \opn\ord{ord}
\opn\Gin{Gin} \opn\Hilb{Hilb}\opn\sort{sort} \opn\Gale{Gale}
\opn\aff{aff} \opn\conv{conv} \opn\relint{relint} \opn\st{st}   \opn\cone{cone}
\opn\lk{lk} \opn\cn{cn} \opn\core{core} \opn\vol{vol}
\opn\link{link} \opn\star{star}\opn\lex{lex} \opn\Gr{Gr}
\opn\gr{gr}

\def\pot#1#2{#1[\kern-0.28ex[#2]\kern-0.28ex]}

\opn\dirlim{\underrightarrow{\lim}}
\opn\inivlim{\underleftarrow{\lim}}
\def\Implies{\ifmmode\Longrightarrow \else
        \unskip${}\Longrightarrow{}$\ignorespaces\fi}
\def\implies{\ifmmode\Rightarrow \else
        \unskip${}\Rightarrow{}$\ignorespaces\fi}
\def\iff{\ifmmode\Longleftrightarrow \else
        \unskip${}\Longleftrightarrow{}$\ignorespaces\fi}

\let\:=\colon
\newtheorem{Theorem}{Theorem}[section]

\newtheorem{Corollary}[Theorem]{Corollary}
\newtheorem{Proposition}[Theorem]{Proposition}
\newtheorem{Remark}[Theorem]{Remark}

\newtheorem{Example}[Theorem]{Example}

\newtheorem{Definition}[Theorem]{Definition}

\let\epsilon\varepsilon
\let\kappa=\varkappa
\makeatletter

\def\qed{\ifhmode\textqed\fi
      \ifmmode\ifinner\quad\qedsymbol\else\dispqed\fi\fi}
\def\textqed{\unskip\nobreak\penalty50
       \hskip2em\hbox{}\nobreak\hfil\qedsymbol
       \parfillskip=0pt \finalhyphendemerits=0}
\def\dispqed{\rlap{\qquad\qedsymbol}}

\opn\dis{dis}
\def\pnt{{\raise0.5mm\hbox{\large\bf.}}}

\opn\Lex{Lex}

\begin{document}


\title{Splittings of toric ideals of graphs}

\author[1]{Anargyros Katsabekis}
\thanks{Corresponding author: Anargyros Katsabekis}

\author[2]{Apostolos Thoma}

\address{Anargyros Katsabekis, Department of Mathematics, University of Ioannina, Ioannina 45110, Greece}
\email{katsampekis@uoi.gr}

\address{Apostolos Thoma, Department of Mathematics, University of Ioannina, Ioannina 45110, Greece}
\email{athoma@uoi.gr}

\keywords{Splittings, toric ideals, graphs, minimal systems of binomial generators.}
\subjclass{13F65, 14M25, 05C25, 05E40.}

\begin{abstract} Let $G$ be a simple graph on the vertex set $\{v_{1},\ldots,v_{n}\}$. An algebraic object attached to $G$ is the toric ideal $I_G$. We say that $I_G$ is subgraph splittable if there exist subgraphs $G_1$ and $G_2$ of $G$ such that $I_G=I_{G_1}+I_{G_2}$, where both $I_{G_1}$ and $I_{G_2}$ are not equal to $I_G$. We show that $I_G$ is subgraph splittable if and only if it is edge splittable. We also prove that the toric ideal of a complete bipartite graph is not subgraph splittable. In contrast, we show that the toric ideal of a complete graph $K_n$ is always subgraph splittable when $n \geq 4$. Additionally, we show that the toric ideal of $K_n$ has a minimal splitting if and only if $4 \leq n \leq 5$. Finally, we prove that any minimal splitting of $I_G$ is also a reduced splitting.
\end{abstract}
\maketitle

\section{Introduction}

Let $A=\{{\bf a}_{1},\ldots,{\bf a}_{m}\}$ be a subset of $\mathbb{N}^{n}$. Consider the polynomial ring $K[x_{1},\ldots,x_{m}]$ over any field $K$. We grade $K[x_{1},\ldots,x_{m}]$ by the semigroup $\mathbb{N}A=\{l_{1}{\bf a}_{1}+\cdots+l_{m}{\bf a}_{m} \mid l_{1},\ldots,l_{m} \in \mathbb{N}\}$ setting ${\rm deg}_{A}(x_{i})={\bf a}_{i}$ for $i=1,\ldots,m$. The $A$-{\em degree} of a monomial ${\bf x}^{\bf u}=x_{1}^{u_1} \cdots x_{m}^{u_m}$ is defined by ${\rm deg}_{A}({\bf x}^{\bf u})=u_{1}{\bf a}_{1}+\cdots+u_{m}{\bf a}_{m} \in \mathbb{N}A$. The {\em $A$-fiber} of a vector ${\bf b} \in \mathbb{N}^{n}$ is the set of all monomials in $K[x_{1},\ldots,x_{m}]$ with $A$-degree equal to ${\bf b}$. The {\em toric ideal} $I_A$ is the prime ideal generated by all the binomials ${\bf x}^{\bf u}-{\bf x}^{\bf v}$ such that ${\rm deg}_{A}({\bf x}^{\bf u})={\rm deg}_{A}({\bf x}^{\bf v})$. 

A recent direction in the theory of toric ideals is to decide when $I_A$ is splittable, see \cite{FHKT, GS, K}. The toric ideal $I_A$ is {\em splittable} if it has a toric splitting, i.e., if there exist toric ideals $I_{A_1}, I_{A_2}$ such that $I_A=I_{A_1}+I_{A_2}$ and $I_{A_i}\not =I_{A}$, for all $1 \leq i \leq 2$. To every simple graph $G$ one can associate the toric ideal $I_G$. In \cite{FHKT} G. Favacchio, J. Hofscheier, G. Keiper, and A. Van Tuyl consider the aforementioned problem for the toric ideal $I_G$, namely they study when $I_G$ has a toric splitting of the form $I_G=I_{G_1}+I_{G_2}$, where $I_{G_i}$, $1 \leq i \leq 2$, is the toric ideal of a subgraph $G_i$ of $G$. More precisely, given a graph $G$ and an even cycle $C$, they consider the graph $H$ which is formed by identifying any edge of $G$ with an edge of $C$. They show \cite[Theorem 3.7]{FHKT} that $I_{H}=I_{G}+I_{C}$ is a splitting of $I_H$. They also prove \cite[Corollary 4.8]{FHKT} that if $G_1$ and $G_2$ form a splitting of $G$ along an edge $e$ and at least one of $G_1$ or $G_2$ is bipartite, then $I_{G}=I_{G_1}+I_{G_2}$ is a splitting of $I_{G}$. Moreover, when $I_G$ has such a splitting, they show \cite[Theorem 4.11]{FHKT} how the graded Betti numbers of $I_G$ are related with those of $I_{G_1}$ and $I_{G_2}$. P. Gimenez and H. Srinivasan showed \cite[Theorem 3.4]{GS} that if $G_1$ and $G_2$ form a splitting of $G$ along an edge $e$, then $I_G$ splits into $I_{G}=I_{G_1}+I_{G_2}$ if and only if at least one of $G_1$ or $G_2$ is bipartite.

This paper aims to partially answer \cite[Question 5.1]{FHKT}, namely for what graphs $G$ can we find $G_1$ and $G_2$ so that their respective toric ideals satisfy $I_G=I_{G_1}+I_{G_2}$? More generally, can we classify when $I_G$ is a splittable toric ideal in terms of $G$? We give a complete answer to the latter question in the case that both $G_1, G_2$ are subgraphs of $G$. We call such a splitting a subgraph splitting and the ideal $I_G$ is called subgraph splittable. Note that there may be splittings of graphs $G$ which are not subgraph splittings, see Example \ref{Nosubgraphsplittins}. All the splittings discussed in the paper pertain to cases where $G_1$ and $G_2$ are subgraphs of $G$.
Our approach is based on the graphs $G\symbol{92}e$ and $G_{S}^e$ introduced in Section 2, where $e$ is an edge of $G$ and $S$ is a minimal system of binomial generators of $I_G$. We show that $I_{G}$ is subgraph splittable if and only if there is an edge $e$ of $G$ and a minimal generating set of binomials $S$ of $I_G$ such that $I_G=I_{G_S^e}+I_{G\symbol{92}e}$ is a splitting, see Theorem \ref{basicsplit}. As an application of our results, we prove that the toric ideal of a complete bipartite graph is not subgraph splittable (see Corollary \ref{Completebipartite}) and the toric ideal of the wheel graph is subgraph splittable if and only if either $n=4$ or $n$ is odd, see Theorem \ref{Suspension}. In Section 3 we study the case that $G$ coincides with the complete graph $K_n$ on $n$ vertices. We show that $I_{K_n}$ is subgraph splittable if and only if $n \geq 4$, see Theorem \ref{SplittableComplete}. Moreover, we introduce minimal splittings and show (Theorem \ref{MinimalSplitting}) that $I_{K_n}$ does not have a minimal splitting for $n \geq 6$. In Section 4 we define reduced splittings of toric ideals and prove that every minimal splitting of $I_G$ is also reduced, see Theorem \ref{ReducedMinimal}.

\section{Edge splittings}
\label{section 2}

In this section, we first collect important notations and definitions used in the paper. For unexplained terminology in graph theory, we refer to \cite{RTT}. Let $G$ be a finite, connected and undirected graph having no loops and no multiple edges on the vertex set $V(G)=\{v_{1},\ldots,v_{n}\}$, and let $E(G)=\{e_{1},\ldots,e_{m}\}$ be the set of edges of $G$. Two edges of $G$ are called {\em adjacent} if they share a common vertex. To each edge $e=\{v_{i},v_{j}\} \in E(G)$, we associate the vector ${\bf a}_{e} \in \{0,1\}^{n}$ defined as follows: the $i$th entry is $1$, the $j$th entry is 1, and the remaining entries are zero. By $I_G$ we denote the toric ideal $I_{A_{G}}$ in $K[e_{1},\ldots,e_{m}]$, where $A_{G}=\{{\bf a}_{e}| e \in E(G)\} \subset \mathbb{N}^{n}$.

A {\em walk} of length $q$ of $G$ connecting $v_{1} \in V(G)$ with $v_{q+1}$ is a finite sequence of the form
$w=(\{v_{1},v_{2}\}, \{v_{2},v_{3}\},\ldots,\{v_{q-1},v_{q}\}, \{v_{q},v_{q+1}\})$ with each $\{v_{i},v_{i+1}\} \in E(G)$, $1 \leq i \leq q$. An {\em even} (respectively, {\em odd}) walk is a walk of even (respectively, odd) length. The walk $w$ is called {\em closed} if $v_{q+1}=v_{1}$. A {\em cycle} is a closed walk
$$w=(\{v_{1},v_{2}\}, \{v_{2},v_{3}\},\ldots,\{v_{q-1},v_{q}\}, \{v_{q},v_{1}\})$$ with $q \geq 3$ and $v_{i}\not = v_{j}$, for every $1 \leq i<j \leq q$. A {\em cut vertex} is a vertex of $G$ whose removal increases the number of connected components of the remaining subgraph. A graph is called {\em biconnected} if it is connected and does not contain a cut vertex. A {\em block} is a maximal biconnected subgraph of $G$.

Given an even closed walk $w = (e_{i_1}, e_{i_2}, \ldots, e_{i_{2q}})$ of $G$, we write $B_w$ for the binomial $B_{w}=\prod_{k=1}^{q}e_{i_{2k-1}}-\prod_{k=1}^{q}e_{i_{2k}} \in I_{G}$. By \cite[Proposition 10.1.5]{Vil} the ideal $I_G$ is generated by all the binomials $B_{w}$, where $w$ is an even closed walk of $G$. We say that $w$ is a {\em primitive walk} if the corresponding binomial $B_w$ is primitive. Recall that given a set of vectors $A \subset \mathbb{N}^n$, the binomial ${\bf x}^{\bf u}-{\bf x}^{\bf v}$ in $I_A$ is called {\em primitive} if there exists no other binomial ${\bf x}^{\bf w}-{\bf x}^{\bf z} \in I_A$ such that ${\bf x}^{\bf w}$ divides ${\bf x}^{\bf u}$ and ${\bf x}^{\bf z}$ divides ${\bf x}^{\bf v}$. Every minimal binomial generator of $I_A$ is primitive, see \cite{St}. A necessary characterization of the primitive elements of a graph $G$ was given in \cite[Lemma 2.1]{Ohsugi-Hibi}:
If $B\in I_{G}$ is primitive, then we have $B=B_{w}$ where $w$ is
one of the following even closed walks:
\begin{enumerate}
  \item $w$ is an even cycle of $G$
  \item $w=(c_{1},c_{2})$, where $c_{1}$ and $c_{2}$ are odd cycles of $G$ having exactly one common vertex
  \item $w=(c_{1},w_{1},c_{2},w_{2})$, where $c_{1}$ and $c_{2}$ are odd cycles of $G$ having no
  common vertex and where $w_{1}$ and $w_{2}$ are walks of $G$ both of which combine a vertex
  $v_{1}$ of $c_{1}$ and a vertex $v_{2}$ of $c_{2}$.
\end{enumerate}

Every even primitive walk $w = (e_{i_1}, e_{i_2}, \ldots, e_{i_{2q}})$ partitions the set of edges of $w$ into two sets $E_{1}=\{e_{i_j}| j \ \textrm{odd}\}$ and $E_{2}=\{e_{i_j}| j \ \textrm{even}\}$. The edges of $E_{1}$ are called {\em odd edges} of $w$ and those of $E_{2}$ are called {\em even edges}. A {\em sink} of a block $B$ is a common vertex of two odd or two even edges of the primitive walk $w$ which belongs to $B$. The primitive walk $w$ is called {\em strongly primitive} if it has no two sinks with distance one in any cyclic block.

Let $w = (e_{1}, e_{2}, \ldots, e_{2q})$ be an even primitive walk and $f=\{v_{i},v_{j}\} \in E(G) \setminus E(w)$ be a chord of $w$ with $i<j$. Then $f$ breaks $w$ in two walks: $w_{1}=(e_{1},\ldots,e_{i-1},f, e_{j},\ldots,e_{2q})$ and $w_{2}=(e_{i},\ldots,e_{j-1},f)$. The chord $f$ is called {\em bridge} of $w$ if there exist two different blocks $B_i$, $B_j$ of $w$ such that $v_{i} \in B_{i}$ and $v_{j} \in B_{j}$. The chord $f$ is called {\em even} (respectively {\em odd}) if it is not a bridge and breaks $w$ in two even walks (respectively odd). 

Let $f=\{v_{i},v_{j}\}$ be an odd chord of $w$ with $i<j$ and $f'=\{v_{k},v_{l}\}$ be another odd chord of $w$ with $k<l$. We say that the odd chords $f$ and $f'$ {\em cross effectively} in $w$ if $k-i$ is odd and either $i<k<j<l$ or $k<i<l<j$. We call an $F_4$ of the walk $w$ a cycle $(e,f,e',f')$ of length four which consists of two edges $e, e'$ of the walk $w$ both odd or both even, and two odd chords $f$ and $f'$ which cross effectively in $w$. We say that the odd chords $f, f'$ cross {\em strongly effectively} in $w$ if they cross effectively and they do not form an $F_4$ in $w$.

A binomial $B_{w} \in I_G$ is called {\em minimal} if it belongs to a minimal system of binomial generators of $I_G$. Since $I_G$ is homogeneous, the graded version of Nakayama's Lemma implies that every minimal system of generators of $I_G$ has the same cardinality. 

The next theorem provides a characterization for the minimal binomials of $I_G$.

\begin{Theorem} \label{Minimalbasic} (\cite[Theorem 4.13]{RTT}) Let $w$ be an even closed walk of $G$. Then $B_w$ is a minimal binomial if and only if \begin{enumerate}
\item[(1)] $w$ is strongly primitive,
\item[(2)] all the chords of $w$ are odd and there are not two of them which cross strongly effectively and
\item[(3)] no odd chord crosses an $F_{4}$ of the walk $w$.
\end{enumerate}
\end{Theorem}

A binomial $B_{w} \in I_G$ is called {\em indispensable} if every system of binomial generators of $I_G$ contains $B_w$ or $-B_{w}$. The next theorem provides a characterization for the indispensable binomials of $I_G$.

\begin{Theorem} \label{Minimalbasic1} (\cite[Theorem 4.14]{RTT}) Let $w$ be an even closed walk of $G$. Then $B_w$ is an indispensable binomial if and only if $w$ is a strongly primitive walk, all the chords of $w$ are odd and there are not two of them that cross effectively.
\end{Theorem}

It follows from Theorems \ref{Minimalbasic} and \ref{Minimalbasic1} that the binomial $B_w$ is not indispensable due to the existence of $F_4$'s in the walk $w$. Moreover, a minimal binomial $B_w$ of $I_G$ is indispensable if and only if the walk $w$ does not have any $F_4$. 

Note that there may exist a subgraph $H$ of the graph $G$ such that $I_H=I_G$. This can happen when there are edges in $G$ that are not used in any walk $w$ such that $B_w$ is a minimal binomial of $I_G$.

Given an edge $e$ of $G$, we denote by $G\symbol{92}e$ the graph with the same vertex set as $G$ and whose edge set consists of all edges of $G$ except $e$. For any set $F \subset E(G)$ of edges of $G$, we use $G\symbol{92} F$ to denote the subgraph of $G$ containing the same vertices as $G$ but with all the edges of $F$ removed.

Let $S=\{B_{w_1}, B_{w_2}, \ldots, B_{w_r}\}$ be a minimal generating set of $I_G$. Given an edge $e$ of $G$, we define $G_S^e$ to be the subgraph of $G$ on the vertex set $$V(G_S^e)=\bigcup_{1 \leq i \leq r \  \text{and}  \ e\in E(w_i)}V(w_i)$$ with edges $$E(G_S^e)=\bigcup_{1 \leq i \leq r \  \text{and} \  e\in E(w_i)}E(w_i).$$
Thus to form the graph $G_S^e$ one needs to first find all binomials $B_{w_i} \in S$, $1 \leq i \leq r$, such that $e$ is an edge of the walk $w_i$. Then we take all vertices and edges of such walks.

We use the symbol $G_S^e$ to emphasize that the graph depends not only on the edge $e$ but also on the minimal system of binomial generators $S$, see Example \ref{ExampleBasic11}. For a toric ideal $I_{G}$ with a unique minimal system of binomial generators $S$, we will simply write $G^{e}$ instead of $G_{S}^{e}$.

\begin{Example} \label{ExampleBasic11} {\rm Recall that the complete graph $K_n$ is the graph with $n$ vertices in which each vertex is connected to every other vertex. Let $G=K_4$ be the complete graph on the vertex set $\{v_{1},\ldots,v_{4}\}$. Let $w_{1}=(\epsilon_{12}, \epsilon_{23}, \epsilon_{34}, \epsilon_{14})$,  $w_{2}=(\epsilon_{12}, \epsilon_{24}, \epsilon_{34}, \epsilon_{13})$ and $w_{3}=(\epsilon_{23}, \epsilon_{13}, \epsilon_{14}, \epsilon_{24})$, where $\epsilon_{ij}=\{v_{i},v_{j}\}$ for $1 \leq i<j \leq 4$. Then $S=\{B_{w_1}=\epsilon_{12}\epsilon_{34}-\epsilon_{23}\epsilon_{14}, B_{w_2}=\epsilon_{12}\epsilon_{34}-\epsilon_{24}\epsilon_{13}\}$ and  $T=\{B_{w_1}, B_{w_3}=\epsilon_{23}\epsilon_{14}-\epsilon_{13}\epsilon_{24}\}$ are minimal generating sets of $I_{G}$. Let $e=\epsilon_{23}$, then $G_{S}^{e}$ is the cycle given by the walk $w_1$ while $G_{T}^{e}$ is the whole graph $G$. }
 \end{Example}
 
The next example shows that, for the same generating set $S$, two different edges $e \neq e'$ can either lead to $G_{S}^{e} \neq G_{S}^{e'}$ or to $G_{S}^{e} = G_{S}^{e'}$.
\begin{figure}[ht]
\begin{center}
\includegraphics[scale=0.15]{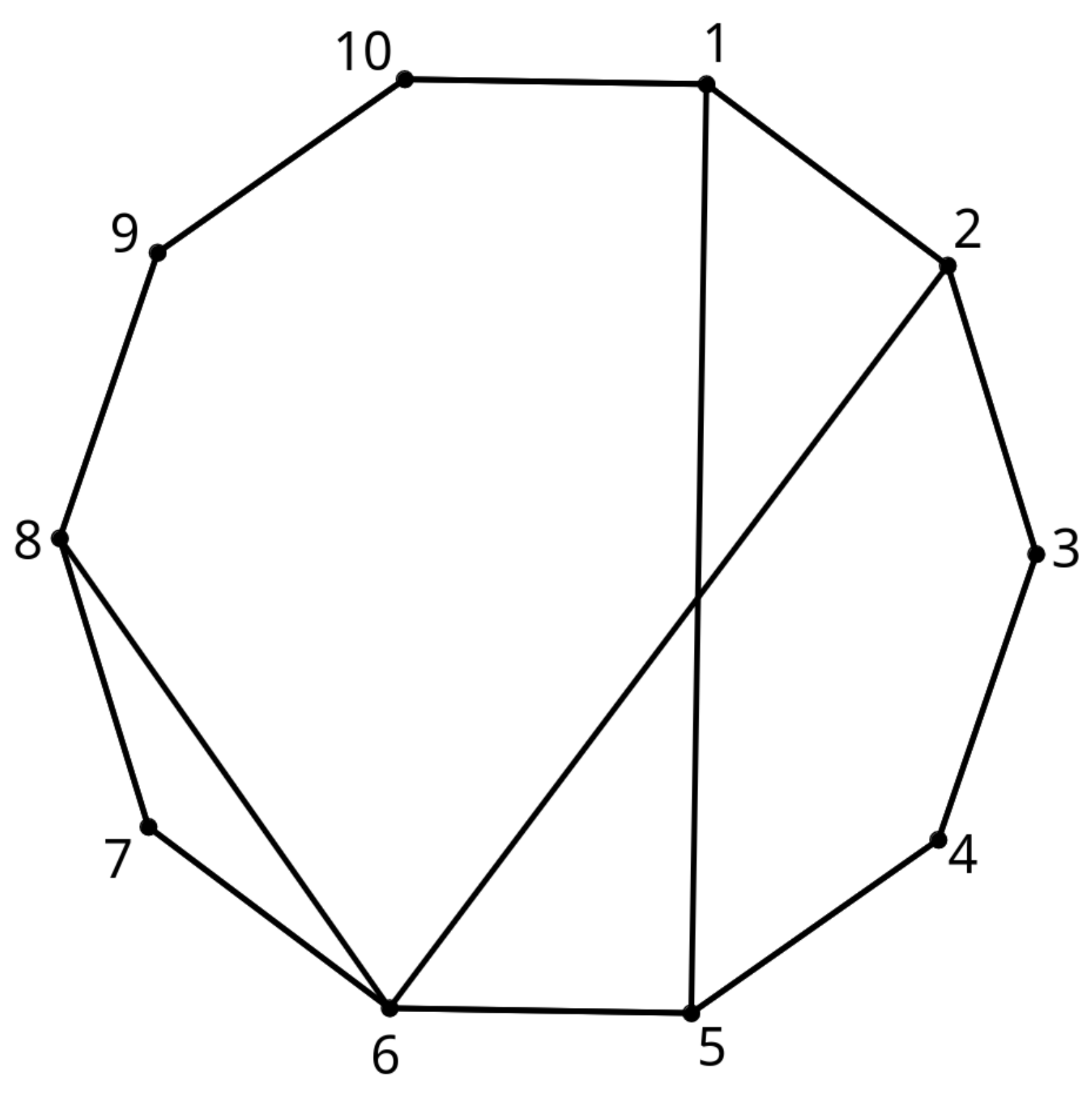}
\caption{The graph $G$ of Example \ref{10gon}.}
\label{10gon}
\end{center}
\end{figure}
 \begin{Example} \label{10gon} {\rm Consider the graph $G$ of Figure \ref{10gon} on the vertex set $\{1,\ldots,10\}$ with edges $e_{1}=\{1,2\}$, $e_{2}=\{2,3\}$, $e_{3}=\{3,4\}$, $e_{4}=\{4,5\}$, $e_{5}=\{5,6\}$, $e_{6}=\{6,7\}$, $e_{7}=\{7,8\}$, $e_{8}=\{8,9\}$, $e_{9}=\{9,10\}$, $e_{10}=\{1,10\}$, $e_{11}=\{1,5\}$, $e_{12}=\{2,6\}$ and $e_{13}=\{6,8\}$. Then $S=\{e_1e_5-e_{11}e_{12}, e_{1}e_{9}e_{13}-e_8e_{10}e_{12}, e_5e_8e_{10}-e_9e_{11}e_{13}, e_2e_4e_6e_{13}-e_{3}e_{5}e_{7}e_{12}, e_2e_4e_6e_8e_{10}-e_{3}e_7e_9e_{11}e_{12}\}$ is a minimal generating set of $I_{G}$. The edge $e_1$ is an edge of the walks that correspond to the binomials $e_1e_5-e_{11}e_{12}, e_{1}e_{9}e_{13}-e_8e_{10}e_{12}$. 
 Thus, $E(G_{S}^{e_1})=\{e_{1}, e_{5}, e_{8}, e_{9}, e_{10}, e_{11}, e_{12}, e_{13}\}$, and $I_{G_{S}^{e_1}}=\langle e_1e_5-e_{11}e_{12}, e_{1}e_{9}e_{13}-e_8e_{10}e_{12}, e_5e_8e_{10}-e_9e_{11}e_{13} \rangle$. Note that $e_1$ is not a variable in the binomial $e_5e_8e_{10}-e_9e_{11}e_{13}$, but this binomial belongs to the ideal $I_{G_{S}^{e_1}}$. Working similarly to the case where $e=e_1$, one can observe that, for $e\in \{e_5, e_8, e_9, e_{10}, e_{11}, e_{12}, e_{13}\}$, we have $G_{S}^{e}=G$. While, for $e\in \{e_2, e_3, e_4, e_6, e_7\}$, we have $G_{S}^{e}=G\symbol{92}{e_1}$.}
 
 \end{Example}
 \begin{Theorem} \label{sum1}
 Let $G$ be a graph and $e$ be an edge of $G$. Let $S=\{B_{w_1}, B_{w_2}, \ldots, B_{w_r}\}$ be a minimal generating set of $I_G$, then $I_G=I_{G_S^e}+I_{G\symbol{92}e}$.
     
 \end{Theorem} 
 {\em \noindent Proof.} We have that $G\symbol{92}e\subset G$ and $G_S^e\subset G$, so $I_{G\symbol{92}e} \subset I_G$ and $I_{G_S^e}\subset I_G$. Thus $I_{G_S^e}+I_{G\symbol{92}e}\subset I_G$. Let $B_{w_i}$, $1 \leq i \leq r$, be a minimal binomial in $S$. Then there are two cases: \begin{enumerate} \item $e\in E(w_i)$. Then $B_{w_i}$ belongs to the ideal $I_{G_S^e}$.
 \item $e \not \in E(w_i)$. Then $B_{w_i}$ belongs to the ideal $I_{G\symbol{92}e}$. 
 \end{enumerate}
 Consequently $I_G\subset I_{G_S^e}+I_{G\symbol{92}e}$. \hfill $\square$
\begin{Example} \label{K_4} {\rm We return to Example \ref{ExampleBasic11}. We have that $I_{G_{S}^{e}}= \langle B_{w_1} \rangle$ and $I_{G_{T}^e}=I_{G}$. Also $G \symbol{92}e$ is the graph with edges $\epsilon_{12}, \epsilon_{34}, \epsilon_{14},  \epsilon_{13}, \epsilon_{24} $, thus $I_{G \symbol{92} e}=\langle B_{w_{2}} \rangle$. Then $I_{G}=I_{G_{S}^{e}}+I_{G \symbol{92}e}$ is a splitting of $I_{G}$, while 
 $I_{G}=I_{G_{T}^e}+I_{G \symbol{92}e}$ is not a splitting of $I_{G}$ since $I_{G_{T}^e}=I_{G}$. }
 \end{Example}
 
\begin{Definition} A subgraph splitting $I_{G}=I_{G_1}+I_{G_2}$ of $I_G$ is called {\em edge splitting} if there exist an edge $e$ of $G$ and a minimal generating set $S$ of $I_G$ such that $G_{1}=G_S^e$ and $G_{2}=G\symbol{92}e$ or $G_{1}=G\symbol{92}e$ and $G_{2}=G_S^e$. The toric ideal $I_G$ is called {\em edge splittable} if there exists an edge splitting of $I_G$.
\end{Definition}

 \begin{Remark} \label{Vertexsplitting} {\rm It is not necessary for all subgraph splittings (if any) of the toric ideal of a graph to be edge splittings. Consider the complete graph $K_5$ on the vertex set $\{v_{1},\ldots,v_{5}\}$. Let $\epsilon_{ij}=\{v_{i},v_{j}\}$ for $1 \leq i<j \leq 5$. Consider the subgraphs $G_1=K_{5} \symbol{92} \{\epsilon_{12}, \epsilon_{34}\}$ and $G_2=K_{5} \symbol{92} \{\epsilon_{14}, \epsilon_{23}\}$ of $K_5$. We have that $S=\{\epsilon_{13}\epsilon_{24}-\epsilon_{14}\epsilon_{23}, \epsilon_{14}\epsilon_{25}-\epsilon_{15}\epsilon_{24}, \epsilon_{23}\epsilon_{45}-\epsilon_{24}\epsilon_{35},\epsilon_{13}\epsilon_{25}-\epsilon_{15}\epsilon_{23},\epsilon_{13}\epsilon_{45}-\epsilon_{14}\epsilon_{35}\}$ is a generating set of $I_{G_1}$ and $T=\{\epsilon_{12}\epsilon_{34}-\epsilon_{13}\epsilon_{24}, \epsilon_{24}\epsilon_{35}-\epsilon_{25}\epsilon_{34}, \epsilon_{13}\epsilon_{45}-\epsilon_{15}\epsilon_{34},\epsilon_{12}\epsilon_{45}-\epsilon_{15}\epsilon_{24}, \epsilon_{12}\epsilon_{35}-\epsilon_{13}\epsilon_{25}\}$ is a generating set of $I_{G_2}$. Also, $S \cup T$ is a generating set of $I_{K_5}$ and $I_{G_i} \neq I_{K_5}$ for all $1 \leq i \leq 2$, so $I_{K_5}=I_{G_1}+I_{G_2}$ is a splitting of $I_{K_5}$ which is not an edge splitting. Note that both $G_1, G_2$ cannot be written in the form $K_{5} \symbol{92} e$.} 
 \end{Remark}

 The next two theorems provide a necessary and sufficient condition for a toric ideal $I_G$ to be subgraph splittable in terms of the graph $G$.
 
 \begin{Theorem} \label{basicsplit1} 
The ideal $I_G$ is edge splittable if and only if there is a minimal system of binomial generators $S=\{B_{w_1},\ldots,B_{w_r}\}$ of $I_G$ with $r \geq 2$ and an edge $e\in E(w_i)$, $1\leq i\leq r$, such that $I_{G_S^e}\not =I_G$. 
     
 \end{Theorem}
 {\em \noindent Proof.} ($\Longrightarrow$) Suppose that the ideal $I_G$ is edge splittable. Then there exist a minimal system of binomial generators $S=\{B_{w_1},\ldots,B_{w_r}\}$ of $I_G$ with $r \geq 2$ and an edge $e$ of $G$ such that $I_G=I_{G_S^e}+I_{G\symbol{92}e}$ is a splitting.
 Thus $I_{G\symbol{92}e}\not =I_G$, so there is $B_{w_i} \in I_G$, $1 \leq i \leq r$, such that $B_{w_i} \notin I_{G\symbol{92}e}$ and therefore $e\in E(w_i)$.\\
  ($\Longleftarrow$) Suppose that there is a minimal system of binomial generators $S=\{B_{w_1},\ldots,B_{w_r}\}$ of $I_G$ with $r \geq 2$ and an edge $e\in E(w_i)$, $1\leq i\leq r$, such that $I_{G_S^e}\not =I_G$. Then from Theorem \ref{sum1} we have that $I_G=I_{G_{S}^e}+I_{G\symbol{92}e}$.
    Since $e\in E(w_i)$, we have that 
     $B_{w_i}\not \in I_{G\symbol{92}e}$ and therefore $I_{G\symbol{92}e}\not =I_G.$ From the hypothesis, it holds that $I_{G_{S}^e}\not =I_G$. Consequently, $I_G$ is edge splittable. \hfill $\square$
                                                   
 Next we state and prove the main result of this article, namely if the toric ideal of a graph has a subgraph splitting then it has also an edge splitting. 

\begin{Theorem} \label{basicsplit} 
The toric ideal $I_G$ is subgraph splittable if and only if it is edge splittable.
     
 \end{Theorem}
{\em \noindent Proof.} ($\Longleftarrow$) If $I_G$ is edge splittable, then $I_G=I_{G_S^e}+I_{G\symbol{92}e}$ is a splitting of $I_G$, and therefore it is subgraph splittable.\\
($\Longrightarrow$) Suppose that $I_G$ is subgraph splittable and let $I_G=I_{G_1}+I_{G_2}$ be a subgraph splitting of $I_G$. Notice that $I_{G_1} \subsetneqq I_G$ and $I_{G_2} \subsetneqq I_G$. Let $\{f_{1},\ldots,f_{s}\}$ be a binomial generating set of $I_{G_1}$ and $\{g_{1},\ldots,g_{t}\}$ be a binomial generating set of $I_{G_2}$. Then $\{f_{1},\ldots,f_{s}, g_{1},\ldots,g_{t}\}$ is a generating set of $I_G$, therefore it contains a minimal binomial generating set $S=\{B_{w_1}, \ldots,B_{w_r}\}$ of $I_G$, since toric ideals of graphs are homogeneous. Note that $r \geq 2$, since $I_{G_1}$ and $I_{G_2}$ are nonzero ideals. But $S \subset \{f_{1},\ldots,f_{s}, g_{1},\ldots,g_{t}\}$, so each $B_{w_i}$, $1 \leq i \leq r$, belongs to at least one of $I_{G_1}, I_{G_2}$. If for every $1 \leq i \leq r$ it holds that $B_{w_i} \in I_{G_2}$, then $I_{G_2}=I_G$ a contradiction. Thus there exists $1 \leq i \leq r$ such that $B_{w_i} \in I_{G_1}$ and $B_{w_i} \notin I_{G_2}$. So there is $e \in E(w_i)$ such that $e \notin E(G_2)$. For every even closed walk $w_j$ such that $e\in E(w_j)$ we have that $B_{w_j}\in I_{G_1}$, thus $I_{G_{S}^e} \subset I_{G_1}$ and therefore $I_{G_{S}^e} \neq I_G$ since $I_{G_1} \subsetneqq I_G$. By Theorem \ref{basicsplit1} the toric ideal $I_G$ is edge splittable. \hfill $\square$\\

 \begin{Example} {\rm In \cite[Example 3.5]{GS} P. Gimenez and H. Srinivasan provide an example of a graph $G$ obtained by gluing two bow ties $G_1$ and $G_2$ along an edge. Since neither $G_1$ nor $G_2$ is bipartite, $I_{G_1}+I_{G_2}$ is not a splitting of $I_G$ by \cite[Theorem 3.4(1)]{GS}. The ideal $I_G$ has a unique minimal system of generators $S$ consisting of five binomials and $I_{G_S^e}=I_G$, for every edge $e$ of $G$. By Theorem \ref{basicsplit1}, the ideal $I_G$ is not edge splittable, and therefore it is not subgraph splittable by Theorem \ref{basicsplit}. Thus $I_G$ does not have a splitting in the form $I_{H_1}+I_{H_2}$ for any subgraphs $H_1, H_2$ of $G$. It is worth mentioning that if a graph $G$ is a gluing of two arbitrary disjoint connected graphs $G_1$ and $G_2$, bipartite or not, along an edge, then from \cite[Theorem 3.4(3)]{GS} there exists a $3$-uniform hypergraph $H$ such that $I_H=I_{G_1}+I_{G_2}$.}
 \end{Example}

To prove that a toric ideal of a graph $I_G$ is not splittable using Theorem \ref{basicsplit1} one has to prove that for every   minimal system of binomial generators $S=\{B_{w_1},\ldots,B_{w_r}\}$ of $I_G$ and each edge $e\in E(w_i)$, $1\leq i\leq r$, is true that $I_{G_S^e} =I_G$. It is unknown if there are graphs $G$ with at least two 
minimal systems of binomial generators $S, S'$ such that  $I_{G_{S}^{e}}=I_{G}$ for each $e \in E(G)$ and there exist an $e' \in E(G)$ such that $I_{G_{S'}^{e'}}\not=I_{G}$. 
Thus Theorem \ref{basicsplit1} is simpler to apply when $I_G$ has a unique minimal system of binomial generators. In particular, this is true for the toric ideal of a bipartite graph, since it is minimally generated by all binomials of the form $B_w$, where $w$ is an even cycle with no chords, see \cite[Theorem 2.3]{HO}. 

A bipartite graph $G$ is called a {\em complete bipartite} graph if its vertex set can be partitioned into two subsets $V_1$ and $V_2$ such that every edge of $V_1$ is connected to every vertex of $V_2$. It is denoted by $K_{m,n}$, where $m$ and $n$ are the numbers of vertices in $V_1$ and $V_2$ respectively.
The next corollary shows that toric ideals of complete bipartite graphs do not admit a subgraph splitting.

 \begin{Corollary} \label{Completebipartite} The toric ideal of $K_{m,n}$ is not subgraph splittable. 
 \end{Corollary}
 {\em \noindent Proof.} According to Theorem \ref{basicsplit}, it suffices to demonstrate that $I_{K_{m,n}}$ is not edge splittable. Let $V_{1}=\{x_{1},\ldots,x_{m}\}$ and $V_2=\{y_{1},\ldots,y_{n}\}$ be the bipartition of the complete bipartite graph $K_{m,n}$ and $E(K_{m,n})=\{b_{ij}| 1 \leq i \leq m, 1\leq j \leq n\}$, where $b_{ij}=\{x_{i},y_{j}\}$. Then $I_{K_{m,n}}$ is minimally generated by the $2 \times 2$ minors of the matrix $M=(b_{ij})$, see \cite[Proposition 10.6.2]{Vil}. Thus $I_{K_{m,n}}$ is minimally generated by the set $S$ of all binomials $b_{ij}b_{kl}-b_{il}b_{kj}$ which are in the form $B_w$, where $w$ is a cycle in $K_{m,n}$ of length $4$. Since $K_{m,n}$ is bipartite, the set $S$ is the unique minimal system of binomial generators of $I_{K_{m,n}}$.
 
 Notice that if $m=1$ or $n=1$, then $I_{K_{m,n}}=\{0\}$. Moreover if $m=2$ and $n=2$, then $I_{K_{m,n}}$ is minimally generated by $b_{11}b_{22}-b_{12}b_{21}$ and therefore it is not subgraph splittable. Assume that $m \geq 2$, $n \geq 2$ and $I_{K_{m,n}} \neq I_{K_{2,2}}$. Let $e$ be any edge of $K_{m,n}$. We claim that $K_{m,n}^e=K_{m,n}$ which implies the equality $I_{K_{m,n}^e}=I_{K_{m,n}}$. Without loss of generality, we can assume that $e=b_{11}$. Then all the generators of $S$ involving $b_{11}$ are $b_{11}b_{ij}-b_{1j}b_{i1}$, for every $1 < i \leq m$ and $1 < j \leq n$. 
 Thus $b_{11}$, $b_{1j}$, $b_{i1}$,  $b_{ij}$, for every $1 < i \leq m$ and $1 < j \leq n$ belong to $E(K_{m,n}^e)$ so $K_{m,n}^e=K_{m,n}$. From Theorem \ref{basicsplit1} it follows that $I_{K_{m,n}}$ is not edge splittable, and therefore $I_{K_{m,n}}$ is not subgraph splittable by Theorem \ref{basicsplit}. \hfill $\square$\\

Let $G=W_{n+1}$ be the wheel graph obtained by connecting a single universal vertex $v_{n+1}$ to all vertices of the cycle  $C_{n}=(\{v_{1},v_{2}\}, \{v_{2},v_{3}\},\ldots, \{v_{n-1},v_{n}\}, \{v_{1},v_{n}\})$, which has length $n \geq 3$. Thus, $G=W_{n+1}$ 
has $n+1$ vertices and $2n$ edges. The cycle edges are the edges of the cycle $C_n$, and the spokes are the edges $\{v_{i},v_{n+1}\}$ for every $1 \leq i \leq n$.
\begin{Example} \label{W4} {\rm Consider the wheel graph $W_{5}$ with five vertices. Let $\epsilon_{12}=\{v_{1},v_{2}\}, \epsilon_{23}=\{v_{2},v_{3}\}, \epsilon_{34}=\{v_{3},v_{4}\}, \epsilon_{14}=\{v_{1},v_{4}\}$, and $\epsilon_{i5}=\{v_{i},v_{5}\}$ for $1 \leq i \leq 4$, be the edges of $G=W_{5}$. Then $$S=\{\epsilon_{12}\epsilon_{45}-\epsilon_{14}\epsilon_{25}, \epsilon_{12}\epsilon_{35}-\epsilon_{23}\epsilon_{15}, \epsilon_{12}\epsilon_{34}-\epsilon_{23}\epsilon_{14}, \epsilon_{23}\epsilon_{45}-\epsilon_{34}\epsilon_{25}, \epsilon_{15}\epsilon_{34}-\epsilon_{14}\epsilon_{35}\}$$ is a minimal generating set of $I_{G}$. Let $e=\epsilon_{15}$, then  $E(G^{e})=\{\epsilon_{12}, \epsilon_{14}, \epsilon_{15}, \epsilon_{23}, \epsilon_{34}, \epsilon_{35}\}.$ Notice that $\epsilon_{12}\epsilon_{45}-\epsilon_{14}\epsilon_{25} \in I_G$ and $\epsilon_{12}\epsilon_{45}-\epsilon_{14}\epsilon_{25} \notin I_{G^e}$, so $I_{G^e} \neq I_{G}$. Thus $I_{G}=I_{G^e}+I_{G \symbol{92} e}$ is a splitting of $I_{G}$ by Theorem \ref{basicsplit1}, and therefore the graph $W_{5}$ is subgraph splittable.}
    
\end{Example}

The next theorem determines when the toric ideal $I_{W_{n+1}}$ is subgraph splittable, where $W_{n+1}$ is the wheel graph with $n+1$ vertices for $n \geq 3$.

 \begin{Theorem} \label{Suspension} Let $G=W_{n+1}$ be the wheel graph with $n+1$ vertices, where $n \geq 3$. \begin{enumerate} \item Suppose that $n$ is even. Then $I_{W_{n+1}}$ is subgraph splittable if and only if $n=4$.
 \item If $n$ is odd, then $I_{W_{n+1}}$ is subgraph splittable.
 \end{enumerate}
 \end{Theorem}

{\em \noindent Proof.} 
(1) Suppose that $n$ is even. If $n=4$, then $W_5$ is subgraph splittable by Example \ref{W4}. According to Theorem \ref{basicsplit}, it suffices to demonstrate that $I_{W_{n+1}}$ is not edge splittable. Since $n>4$ is even, every odd cycle of $W_{n+1}$ contains the universal vertex, thus, from \cite[Proposition 5.5]{OH1}, there is a bipartite graph $H$ such that $I_{G}=I_H$. The graph $H$ is bipartite, so $I_{H}$ has a unique minimal system of generators and the generators are in the form $B_w$ for even walks that correspond to cycles with no chords. 
 The graph $H$ has vertices $v_1, v_2, \dots, v_{n+1}, v_{n+2}$ and edges the cycle edges and the edges $\{v_i, v_{n+1}\}$ for $i$ odd and $\{v_i, v_{n+2}\}$ for $i$ even, see the proof of \cite[Proposition 5.5]{OH1}. We will prove that for any $e\in E(H)$ we have $H^e=H$. There are three cases to consider, namely (i) $e$ is a cycle edge,
 (ii) $e$ is an edge of the form $\{v_i, v_{n+1}\}$ for $i$ odd, and (iii) $e$ is an edge of the form $\{v_i, v_{n+2}\}$ for $i$ even. It is not restrictive to take $e=\{v_1, v_2\}$ for (i),
$e=\{v_1, v_{n+1}\}$ for (ii), and $e=\{v_2, v_{n+2}\}$ for (iii).
   \begin{enumerate} \item[(i)] $e=\{v_1, v_2\}$. We consider the following even walks $w$
   of length $n$, $4$ or $6$ with no chords that have $e$ as an edge: 
   $$(e=\{v_{1},v_{2}\}, \{v_{2},v_{3}\}, \dots , \{v_{n-1},v_{n}\}, \{v_{n},v_{1}\}),$$
   $$(e=\{v_{1},v_{2}\}, \{v_{2},v_{3}\}, \{v_{3},v_{n+1}\}, \{v_{n+1},v_{1}\}),$$
   $$(e=\{v_{1},v_{2}\}, \{v_{2},v_{n+2}\}, \{v_{n+2},v_{n}\}, \{v_{n},v_{1}\}),$$
   $$(e=\{v_{1},v_{2}\}, \{v_{2},v_{n+2}\}, \{v_{n+2},v_{i}\}, \{v_{i},v_{i+1}\}, \{v_{i+1},v_{n+1}\}, \{v_{n+1},v_{1}\}),$$ where $i$ is even and $4\leq i\leq n-2$. Note that all these walks contain $e$, every edge of $H$ is in at least one of these walks and these walks $w$ do not have chords, therefore the corresponding $B_{w}$'s are minimal binomials of $I_{H}$, and thus $H^{e}=H$.  
  \item[(ii)] $e=\{v_{1},v_{n+1}\}$.  We consider the following even walks $w$
   of length 4 or 6 with no chords that have $e$ as an edge: 
   $(e=\{v_{1},v_{n+1}\}, \{v_{n+1},v_{3}\}, \{v_{3},v_{2}\}, \{v_{2},v_{1}\}),$
   $$(e=\{v_{1},v_{n+1}\}, \{v_{n+1},v_{n-1}\}, \{v_{n-1},v_{n}\}, \{v_{n},v_{1}\}),$$
   $$(\{v_{1},v_{2}\}, \{v_{2},v_{n+2}\}, \{v_{n+2},v_{i}\}, \{v_{i},v_{i+1}\}, \{v_{i+1},v_{n+1}\}, e=\{v_{n+1},v_{1}\}),$$ where $i$ is even and $4\leq i\leq n-2$, 
   $$(e=\{v_{1},v_{n+1}\}, \{v_{n+1},v_{j}\}, \{v_{j},v_{j+1}\}, \{v_{j+1},v_{n+2}\}, \{v_{n+2},v_{n}\}, \{v_{n},v_{1}\}),$$ where $j$ is odd and $3\leq j\leq n-3$. Note that all these walks contain $e$, every edge of $H$ is in at least one of these walks and these walks $w$ do not have chords, therefore the corresponding $B_{w}$'s are minimal binomials of $I_{H}$, and thus $H^{e}=H$. 
   \item[(iii)] $e=\{v_{2},v_{n+2}\}$. We consider the following even walks $w$ of length 4 or 6 with no chords that have $e$ as an edge: $(e=\{v_{2},v_{n+2}\}, \{v_{n+2},v_{4}\}, \{v_{4},v_{3}\}, \{v_{3},v_{2}\}),$
   $$(e=\{v_{2},v_{n+2}\}, \{v_{n+2},v_{n}\}, \{v_{n},v_{1}\}, \{v_{1},v_{2}\}),$$
   $$(\{v_{1},v_{2}\}, e=\{v_{2},v_{n+2}\}, \{v_{n+2},v_{i}\}, \{v_{i},v_{i+1}\}, \{v_{i+1},v_{n+1}\}, \{v_{n+1},v_{1}\}),$$ where $i$ is even and $4\leq i\leq n-2$,
   $$(e=\{v_{n+2},v_{2}\}, \{v_{2},v_{3}\}, \{v_{3},v_{n+1}\}, \{v_{n+1},v_{j}\}, \{v_{j},v_{j+1}\}, \{v_{j+1},v_{n+2}\}),$$ where $j$ is odd and $5\leq j\leq n-1$. Note that all these walks contain $e$, every edge of $H$ is in at least one of these walks and these walks $w$ do not have chords, therefore the corresponding $B_{w}$'s are minimal binomials of $I_{H}$, and thus $H^{e}=H$.
\end{enumerate}
Thus in all cases $H^{e}=H$. The edges of the two graphs $H$ and $G$ are in one-to-one correspondence and their ideals are equal, see \cite[Proposition 5.5]{OH1}, consequently $G=G^{e}$ for every $e\in G$, so from Theorems \ref{basicsplit1} and \ref{basicsplit} the ideal $I_{G}$ is not subgraph splittable.
   
(2) Suppose that $n$ is odd. For $n=3$ we have that $W_{4}$ is the complete graph on the vertex set $\{v_{1},\ldots,v_{4}\}$, so $I_{W_{4}}$ is splittable by Example \ref{K_4}. Suppose that $n \geq 5$. There are several even cycles, all of them with chords. All cycles of length greater than 4 have at least one even chord, a spoke, thus the corresponding $B_w$'s are not minimal binomials by Theorem \ref{Minimalbasic}. If $w$ is a cycle of length $4$, then it is in the form $(\{v_{n+1},v_{i}\}, \{v_{i},v_{i+1}\}, \{v_{i+1},v_{i+2}\}, \{v_{i+2},v_{n+1}\}),$ and has exactly one odd chord, namely $\{v_{n+1},v_{i+1}\}$. Thus $B_w$ is an indispensable binomial of $I_G$ by Theorem \ref{Minimalbasic1}. Notice that two of the edges of $w$ are spokes and the other two are two consecutive cycle edges. 

Any odd cycle of $G$ either coincides with $C_{n}$ or has at least three vertices, namely the universal vertex $v_{n+1}$ and at least two vertices of $C_{n}$. Thus $G$ has no two odd vertex disjoint cycles. 

Let $B_{w} \in I_{G}$ be a minimal binomial, then it is also primitive. Since $G$ has no two odd vertex disjoint cycles by
\cite[Lemma 2.1]{Ohsugi-Hibi}, see also the introduction of section 2, $w$ is an even cycle or $w=(c_1, c_2)$, where $c_{1}$ and $c_{2}$ are odd cycles of $G$ having exactly one common vertex and none of the cycles $c_{1}, c_{2}$ is $C_{n}$. In the case of wheel graphs this common vertex is the universal 
vertex $v_{n+1}$. If one of the cycles has length more than three, then it has at least two chords which are spokes and actually both are bridges of the walk, thus the binomial $B_w$ is not minimal by Theorem \ref{Minimalbasic}.  

Thus we conclude that, for a minimal binomial $B_w$, there are two cases: \begin{enumerate} \item $w$ is a cycle of length $4$ with exactly one odd chord. Then $B_w$ is an indispensable binomial of $I_G$ by Theorem \ref{Minimalbasic1}, and two of the edges of the cycle $w$ are spokes and the other two are two consecutive cycle edges. 
\item $w$ has no chords or bridges and it is of the form $w=(c_{1},c_{2})$ where $c_1$, $c_2$ are odd cycles of length 3 intersecting in exactly one vertex, namely $v_{n+1}$. 
There exist minimal binomial generators like these only for $n>5$. By Theorem \ref{Minimalbasic1} the binomial $B_w$ is an indispensable binomial of $I_G$.
\end{enumerate}
Thus $I_G$ has a unique minimal system of generators.

Let $e=\{v_{1}, v_{2}\}$. We claim that $\epsilon_{34}=\{v_{3},v_{4}\}$ does not belong to $E(G^e)$. The cycle edges $e$ and $\epsilon_{34}$ are not consecutive thus there is no
minimal generator of type (i), $B_w$ where $w$ is a cycle of length four, that 
contains the edges $e$ and $\epsilon_{34}$. In the binomials $B_w$ of type (ii) in the walk
$w=(c_{1},c_{2})$ each odd cycle $c_1$, $c_2$ of length 3 contains exactly one cycle edge.
Thus the only possible walk  that contains the edges $e$ and $\epsilon_{34}$ is
$\gamma=(\gamma_{1},\gamma_{2})$, where $\gamma_{1}=(e=\{v_{1}, v_{2}\},\{v_{2},v_{n+1}\},\{v_{1},v_{n+1}\})$ and $\gamma_{2}=(\epsilon_{34},\{v_{4},v_{n+1}\},\{v_{3},v_{n+1}\})$. But then $B_{\gamma}$ is not a minimal binomial of $I_{G}$ since there is a bridge: $\{v_{2},v_{3}\}$. Note that $v_2$ belongs to the block $\gamma_{1}$ and $v_3$ belongs to the block $\gamma_{2}$ of $\gamma$. 

We conclude that  the edge $\epsilon_{34}$ does not belong to $E(G^e)$. Consider the even cycle $\zeta=(\{v_{2},v_{3}\},\{v_{3},v_{4}\}, \{v_{4},v_{n+1}\}, \{v_{2},v_{n+1}\})$, then $B_{\zeta}$ is a minimal binomial of $I_{G}$ and also $B_{\zeta} \notin I_{G^e}$, since $\epsilon_{34}$ is one of the edges of $\zeta$. Thus $I_{G} \neq I_{G^e}$ and therefore $I_{G}$ is splittable by Theorem \ref{basicsplit1}. \hfill $\square$

\section{  Minimal splittings and the complete graph} \label{complete}

In this section, we study the special case of toric ideals of complete graphs. We prove that the toric ideal of a complete graph with more than three vertices is always subgraph splittable and we provide all possible splittings.  We introduce minimal splittings of toric ideals of graphs and prove that, only for the complete graphs with four or five vertices, the corresponding toric ideals have a minimal splitting. 

Let $n \geq 4$ be an integer and $K_n$ be the complete graph on the vertex set $\{v_{1},\ldots,v_{n}\}$ with edges 
$\{\epsilon_{ij}|1 \leq i<j \leq n \}$, where $\epsilon_{ij}=\{v_{i}, v_{j}\}$. In contrast to toric ideals of (complete) bipartite graphs which always have a unique minimal set of binomial generators, toric ideals of complete graphs have a huge number of different minimal systems of binomial generators.
By \cite[Theorem 3.3]{Vil1} the set $$\mathcal{S}:=\{\epsilon_{ij}\epsilon_{kl}-\epsilon_{il}\epsilon_{jk}, \epsilon_{ik}\epsilon_{jl}-\epsilon_{il}\epsilon_{jk}| 1\leq i<j<k<l \leq n\}$$ is a minimal generating set of $I_{K_n}$. Let $\{{\bf e}_{1},\ldots,{\bf e}_{n}\}$ be the canonical basis of $\mathbb{R}^n$. Since the above set is a minimal generating set of $I_{K_n}$, the only $A_{K_n}$-fibers contributing to the minimal generators are those consisting of all monomials with $A_{K_n}$-degree ${\bf e}_{i}+{\bf e}_{j}+{\bf e}_{k}+{\bf e}_{l}$, where $1\leq i<j<k<l \leq n$. There are $\binom{n}{4}$ such fibers and each one consists of three monomials, namely $\epsilon_{ij}\epsilon_{kl}, \epsilon_{il}\epsilon_{jk}$ and $\epsilon_{ik}\epsilon_{jl}$, which have no common factor other than $1$. Therefore to generate the ideal $I_{K_n}$ we need to take any two of the three binomials $\epsilon_{ij}\epsilon_{kl}-\epsilon_{il}\epsilon_{jk}, \epsilon_{ik}\epsilon_{jl}-\epsilon_{il}\epsilon_{jk},
\epsilon_{ij}\epsilon_{kl}-\epsilon_{ik}\epsilon_{jl}$,
for every $1\leq i<j<k<l \leq n$, see \cite{ChKTh, DSS} for more details. Thus every minimal system of binomial generators of $I_{K_n}$ consists of $2\binom{n}{4}$ binomials. By \cite[Theorem 2.9]{ChKTh} the ideal $I_{K_n}$ has $3^{\binom{n}{4}}$ different minimal systems of binomial generators, which is a huge number even for a small $n$.\\

The next remark will be used in the proof of Proposition \ref{directbasic}.

\begin{Remark} \label{Indispe} {\rm The monomials $\epsilon_{ij}\epsilon_{kl}, \epsilon_{il}\epsilon_{jk}$ and $\epsilon_{ik}\epsilon_{jl}$, where $1 \leq i<j \leq n,$  are indispensable monomials, namely each one is a monomial term of at least one binomial in every minimal system of binomial generators of $I_{K_n}$.}
\end{Remark}

\begin{Theorem} \label{SplittableComplete} The toric ideal of $K_n$ is subgraph splittable if and only if $n\ge 4$.  
 \end{Theorem}
 {\em \noindent Proof.} For $n\in \{1, 2, 3\}$ we have that $I_{K_n}=\{0\}$, so $I_{K_n}$ is not splittable. Suppose that $n\ge 4$ and let $G=K_{n}$ and $e=\epsilon_{12}$. Suppose that $\epsilon_{13}$ is an edge of $G_{\mathcal{S}}^{e}$, so there exists a binomial $B \in \mathcal{S}$ in four variables which contains the variables $\epsilon_{12}$ and $ \epsilon_{13}$. Then $B$ necessarily is in the form $\epsilon_{12}\epsilon_{3i}-\epsilon_{13}\epsilon_{2i}$ for an $i\not \in \{1,2,3\}$. But $1 < 2 < 3 <i \leq n$  therefore the only binomials that belong to $\mathcal{S}$ involving the variables  $\epsilon_{12}, \epsilon_{13}, \epsilon_{2i}, \epsilon_{3i}$ are $\epsilon_{12}\epsilon_{3i}-\epsilon_{1i}\epsilon_{23}$ and $\epsilon_{13}\epsilon_{2i}-\epsilon_{1i}\epsilon_{23}$ in $\mathcal{S}$. Thus $B$ does not belong to $\mathcal{S}$, a contradiction. Consequently $\epsilon_{13}$ is not an edge of $G_{\mathcal{S}}^{e}$.  Since $\epsilon_{12}\epsilon_{34}-\epsilon_{13}\epsilon_{24} \in I_{G}$ and $\epsilon_{12}\epsilon_{34}-\epsilon_{13}\epsilon_{24} \not \in I_{G_{\mathcal{S}}^{e}}$, we get $I_{G_{\mathcal{S}}^{e}} \not =I_G$ and therefore $I_{G}=I_{G_{\mathcal{S}}^{e}}+I_{G\symbol{92} e}$ is a splitting of $I_G$ by Theorem \ref{basicsplit1}. \hfill $\square$ 

\begin{Definition} Let $I_G=I_{G_1}+I_{G_2}$ be a subgraph splitting of $I_G$, $S=\{f_1,\dots, f_r\}$ be a minimal system of binomial generators of $I_{G_1}$ and $T=\{g_1,\dots, g_t\}$ be a minimal system of binomial generators of $I_{G_2}$. We say that the splitting $I_G=I_{G_1}+I_{G_2}$ is a {\em minimal splitting} of $I_G$ if $\{f_1,\dots, f_r, g_1,\dots, g_t\}$ is a minimal system of generators of $I_{G}$.
 \end{Definition}

\begin{Remark} {\em (1) The property of being minimal splitting does not depend on the minimal systems of generators chosen in the definition. Suppose that $\{f'_1,\dots, f'_r\}$ is a minimal system of binomial generators of $I_{G_1}$ and $\{g'_1,\dots, g'_t\}$ is a minimal system of binomial generators of $I_{G_2}$. But $I_G=I_{G_1}+I_{G_2}$, so $\{f'_1,\dots, f'_r, g'_1,\dots, g'_t\}$ is a generating set of $I_G$ consisting of $r+t$ elements and therefore it is a minimal system of generators.\\ 
(2) Let $I_G=I_{G_1}+I_{G_2}$ be a minimal splitting of $I_G$. For any $f_i \in S$, $1 \leq i \leq r$, we have that neither $f_i$ nor $-f_i$ belongs to $T$. For any $g_j \in T$, $1 \leq j \leq t$, we have that neither $g_j$ nor $-g_j$ belongs to $S$.\\
(3) All the splittings which appeared in \cite{FHKT, GS} are minimal splittings.}
\end{Remark}
Example \ref{ExampleG_n} presents a graph and an edge splitting of the corresponding toric ideal which is not minimal.   

 \begin{figure}[ht]
\begin{center}
\includegraphics[scale=0.35]{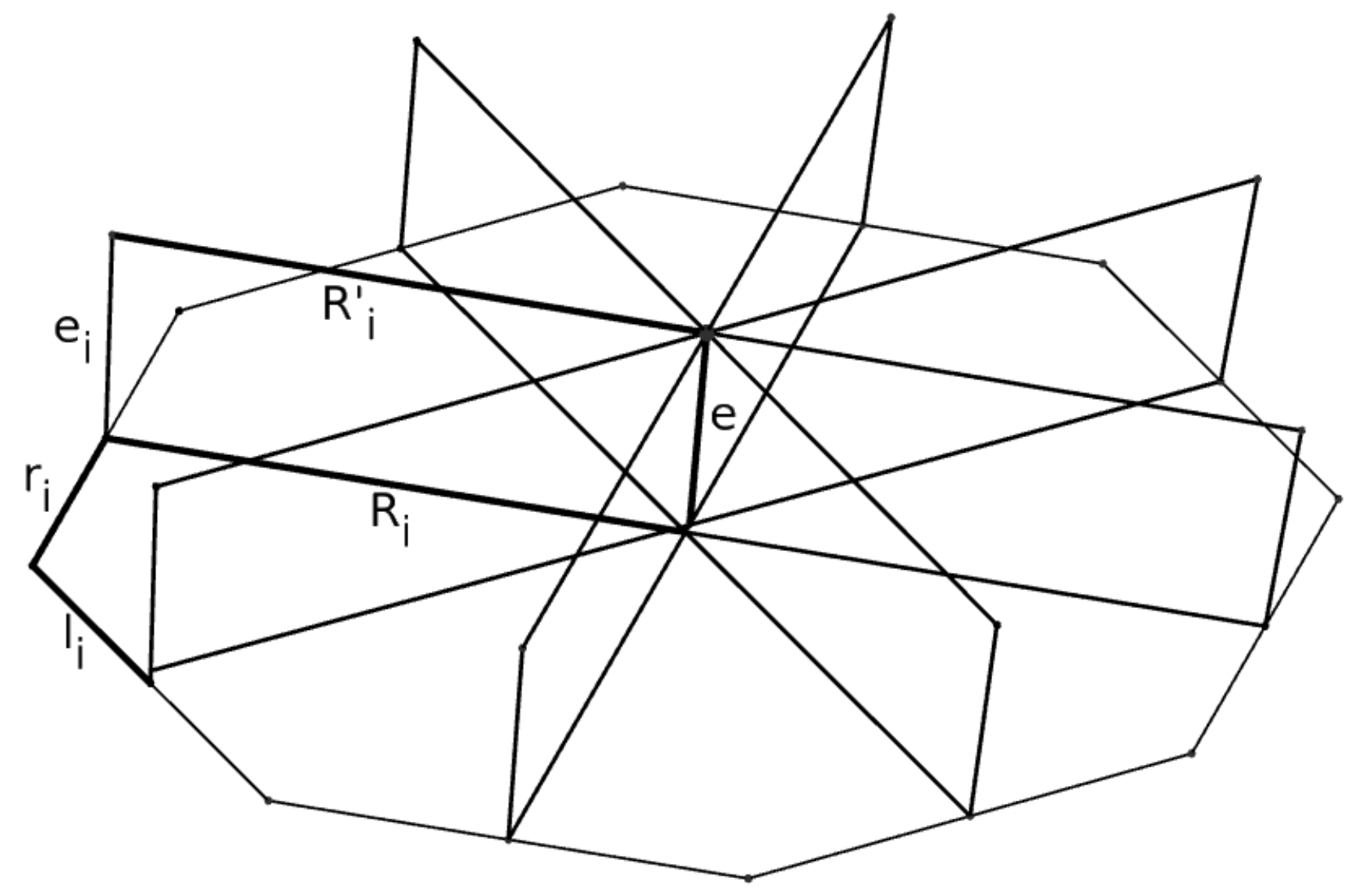}
\caption{The graph $G_8$}
\label{graph}
\end{center}
\end{figure}
 
 \begin{Example} \label{ExampleG_n} {\rm Let $G_n$ be a  graph with $3n+2$ vertices, $$V(G_n)=\{v_1, v_2, \dots, v_n, m_1, m_2, \dots, m_n, m'_1, m'_2, \dots, m'_n, o, o' \},$$ and $5n+1$ edges in the form: $e=\{ o, o'\}$, $e_i=\{m_i, m'_i\}$, $R_i=\{o, m_i\}$, $R'_i=\{o', m'_i\}$, $r_i=\{v_i, m_i\}$, $l_i=\{v_i, m_{i-1}\}$, where $1\leq i \leq n$ and $i-1$ equals $n$ when $i=1$. The roundabout graph $G_n$ is bipartite and also, from \cite[Theorem 2.3]{HO}, we have that $I_{G_n}$ has a unique minimal system of generators consisting of all binomials $B_w$, where $w$ is an even cycle with no chords. The graph $G_n$ has $n$ cycles with no chords of the form $(e, R_i, e_i, R'_i)$, $n$ cycles of the form $(R_i, l_i, r_i, R_{i-1})$, $1$ cycle of the form $(l_1,r_1,l_2, r_2,\ldots,l_n,r_n)$ and $2\binom{n}{2}$ cycles of two forms, namely $(R'_i, e_i, r_i, l_i, r_{i-1}, \ldots, e_j, R'_j)$ and $(R'_i, e_i, l_{i+1}, r_{i+1}, \ldots, e_j, R'_j)$, where $i,j\in \{1,\cdots ,n\}$, $i-1$ equals $n$ when $i=1$ and $i+1$ equals 1 when $i=n$. Thus, $G_n$ has exactly $2n+1+2\binom{n}{2}$ cycles without a chord. 

\begin{figure}[ht]
\begin{center}
\includegraphics[scale=0.35]{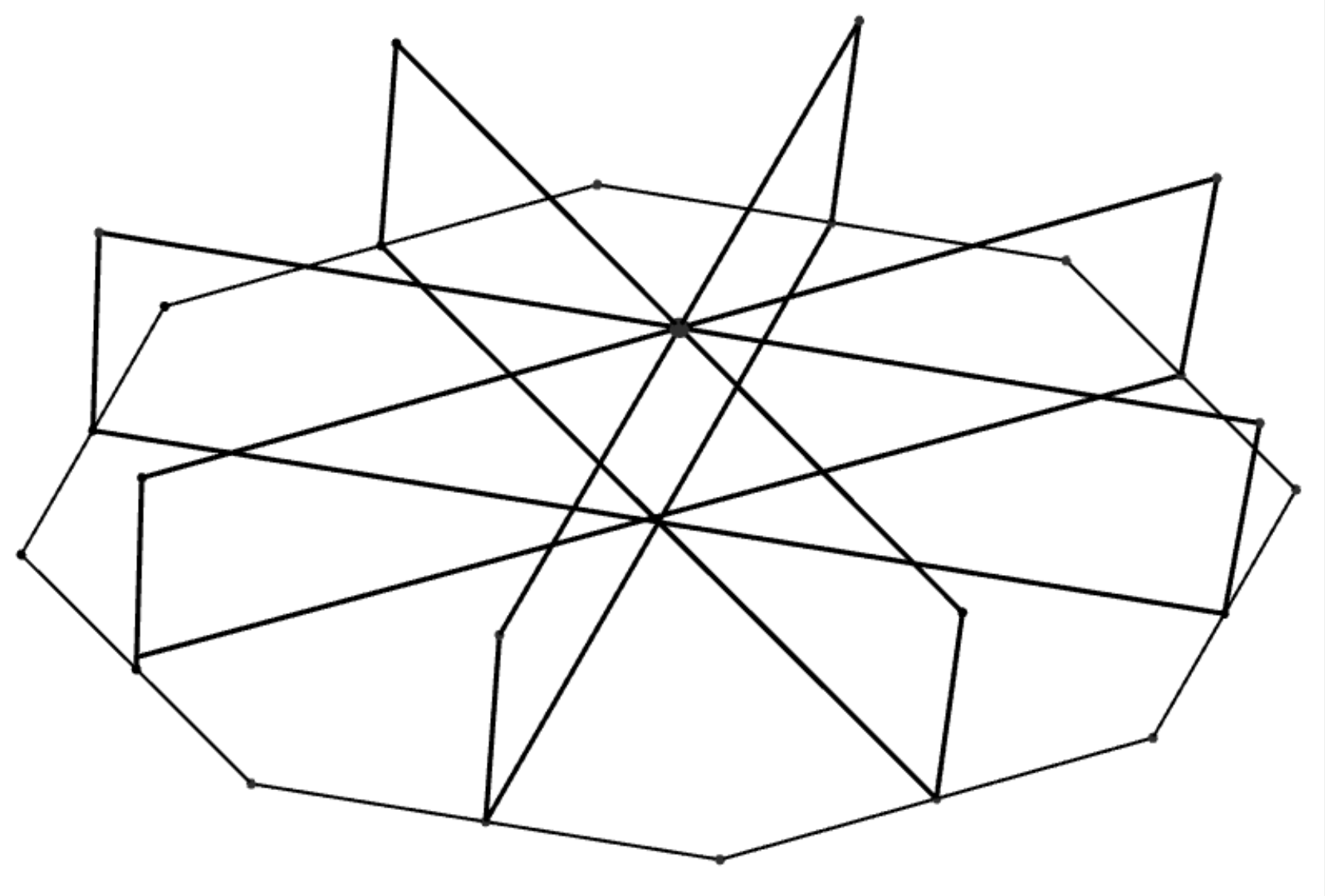}
\caption{The graph ${G_8}\symbol{92}e$}
\label{graph2}
\end{center}
\end{figure}

\begin{figure}[ht]
\begin{center}
\includegraphics[scale=0.35]{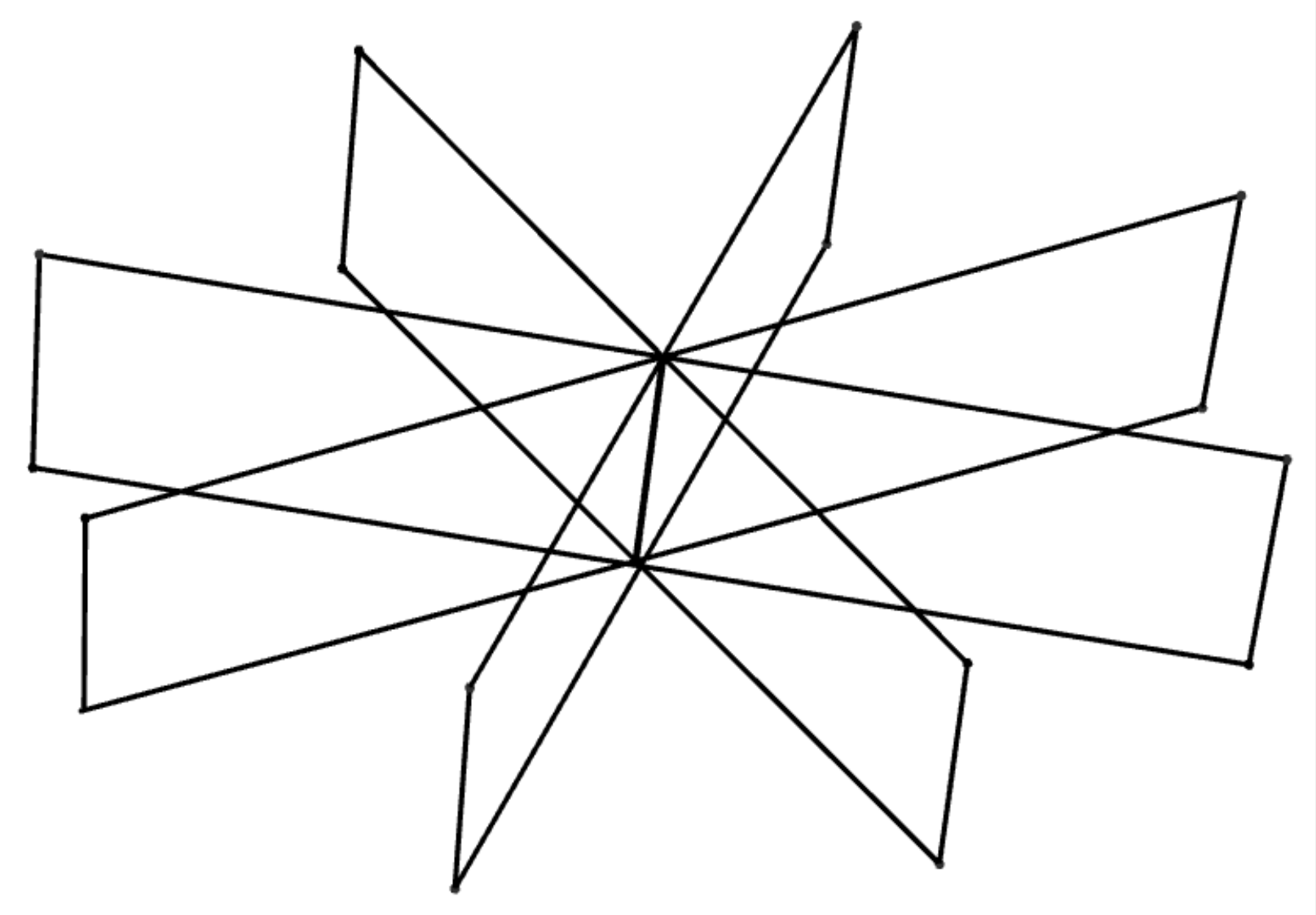}
\caption{The graph $G_8^e$}
\label{graph3}
\end{center}
\end{figure}
 The graph ${G_n}\symbol{92}e$ has $n+3\binom{n}{2}+1$ even cycles with no chords of the following form: $n$ cycles of the form $(R_i, l_i, r_i, R_{i-1})$, $1$ cycle of the form $(l_1, r_1, l_2, r_2,\ldots, l_n, r_n)$, $2\binom{n}{2}$ cycles of two forms, namely $(R'_i, e_i, r_i, l_i, r_{i-1}, \ldots, e_j, R'_j)$ and $(R'_i, e_i, l_{i+1}, r_{i+1}, \ldots, e_j, R'_j)$, and $\binom{n}{2}$ cycles which correspond to cycles $(R_i, e_i, R'_i, R'_j, e_j, R_j)$ in ${G_n}\symbol{92}e$ of length 6 that in the graph $G_n$ had $e$ as an even chord. In general, $I_{{G_n}\symbol{92}e}$ has $n+3\binom{n}{2}+1$ minimal binomials.
 
The graph $G_n^e$ has $n$ cycles with no chords in the form $(e, R_i, e_i, R'_i)$. Thus, the ideal $I_{G_n^e}$ has $n$ minimal binomials. It follows from Theorem \ref{basicsplit1} that $I_{G_n}=I_{G_{n}^{e}}+I_{G \symbol{92}e}$ is an edge splitting. This splitting is not minimal, since the ideal $I_{{G_n}\symbol{92}e}$ is generated minimally by $n+1+3\binom{n}{2}$ binomials, while the ideal $I_{G_n^e}$ is generated minimally by $n$ binomials and the ideal $I_{G_n}$ is generated minimally by $2n+1+2\binom{n}{2}$ binomials.
 
In Figure \ref{graph} we plot the graph $G_8$  with 26 vertices and 41 edges. Note that $G_8$ has exactly $73$ cycles without chords, see Figure \ref{graph}. We plot in Figure \ref{graph2} the graph ${G_8}\symbol{92}e$ and in Figure \ref{graph3} the graph $G_8^e$. Note that the toric ideal of the graph ${G_8}\symbol{92}e$ has 93 minimal binomials, while the toric ideal of $G_8^e$ has 8 minimal binomials.}

 \end{Example}              
\begin{Example} \label{MinimalK_4} {\rm Let $K_4$ be the complete graph on the vertex set $\{v_{1},\ldots,v_{4}\}$. Let $G_1=K_{4} \symbol{92} \{\epsilon_{12}, \epsilon_{34}\}$ and $G_2=K_{4} \symbol{92} \{\epsilon_{14}, \epsilon_{23}\}$ be subgraphs of $K_4$, then $I_{K_4}=I_{G_1}+I_{G_2}$ is a minimal splitting of $I_{K_4}$, since $I_{G_1}= \langle \epsilon_{13}\epsilon_{24}-\epsilon_{14}\epsilon_{23} \rangle$, $I_{G_2}=\langle \epsilon_{13}\epsilon_{24}-\epsilon_{12}\epsilon_{34} \rangle$ and
 $I_{K_4}=\langle \epsilon_{13}\epsilon_{24}-\epsilon_{14}\epsilon_{23}, \epsilon_{13}\epsilon_{24}-\epsilon_{12}\epsilon_{34} \rangle$.}
 \end{Example}

 In Theorem \ref{SplittableComplete} we have seen that the ideal $I_{K_n}$ is splittable for $n\geq 4$. The aim of Propositions \ref{redsplitting}, \ref{directbasic}, \ref{all splittings} is to find all possible splittings of $I_{K_n}$. Finally, Theorem \ref{MinimalSplitting} shows that the toric ideal of $K_n$ has no minimal splitting for $n \geq 6$.

\begin{Proposition}\label{redsplitting}
    Let $n \geq 4$ be an integer and $w=(a,b,c,d)$ be an even cycle of length 4 of $K_n$, where $a, b, c, d \in E(G)$. Then $$I_{K_n}=I_{{K_n}\symbol{92} \{a,c\}}+I_{{K_n}\symbol{92} \{b,d\}}$$ is a splitting of $I_{K_n}.$
\end{Proposition}
{\em \noindent Proof.} Since $K_n$ is the complete graph, the cycle $w$ has two chords, namely $e$ and $f$. Then $ac-ef \in I_{K_n}$ does not belong to $I_{{K_n}\symbol{92} \{a,c\}}$, thus $I_{{K_n}\symbol{92} \{a,c\}} \not =I_{K_n}$. Also $I_{{K_n}\symbol{92} \{b,d\}} \not =I_{K_n}$ since $ef-bd \in I_{K_n}$ does not belong to $I_{{K_n}\symbol{92} \{b,d\}}$.

It remains to show that $I_{K_n} \subset I_{{K_n}\symbol{92} \{a,c\}}+I_{{K_n}\symbol{92} \{b,d\}}$. Let $B_{\gamma}$ be a binomial belonging to the set $S$ defined in the proof of Theorem \ref{SplittableComplete}. Clearly if $B_{\gamma}$ belongs to $I_{{K_n}\symbol{92} \{a,c\}}$ or $ I_{{K_n}\symbol{92} \{b,d\}}$, then $B_{\gamma}$ belongs to $I_{{K_n}\symbol{92} \{a,c\}}+I_{{K_n}\symbol{92} \{b,d\}}$. Suppose that $B_{\gamma}$ does not belong to  $I_{{K_n}\symbol{92} \{a,c\}}$ or $ I_{{K_n}\symbol{92} \{b,d\}}$. Then $\gamma$ contains at least one edge from the set $ \{a,c\}$, say $a$, and at least one edge from the set $ \{b,d\}$, say $b$. Let $\gamma=(a,b,c',d')$ and $e', f'$ be the chords of the walk $\gamma$. Notice that $ac'-e'f'$ belongs to $I_{{K_n}\symbol{92} \{b,d\}}$ and $bd'-e'f'$ belongs to $I_{{K_n}\symbol{92} \{a,c\}}$. Thus $$B_{\gamma}=ac'-bd'=(ac'-e'f')-(bd'-e'f')\in I_{{K_n}\symbol{92} \{a,c\}}+I_{{K_n}\symbol{92} \{b,d\}}.$$
So $I_{K_n}=I_{{K_n}\symbol{92} \{a,c\}}+I_{{K_n}\symbol{92} \{b,d\}}$ is a splitting. \hfill $\square$

\begin{Remark} {\rm Let $I_A=I_{A_1}+I_{A_2}$ be a splitting of $I_A$. If there exists a set $A'_1$ such that $I_{A_1}\subset I_{A'_1} \subsetneqq I_A$ then $I_A=I_{A'_1}+I_{A_2}$ is also a splitting. By Proposition \ref{redsplitting}, $I_{K_n}=I_{{K_n}\symbol{92} \{a,c\}}+I_{{K_n}\symbol{92} \{b,d\}}$ is a splitting of $I_{K_n}$, and therefore $I_{K_n}=I_{{K_n}\symbol{92} a}+ I_{{K_n}\symbol{92} b}$, $I_{K_n}=I_{{K_n}\symbol{92} \{a,c\}}+ I_{{K_n}\symbol{92} b}$ and $I_{K_n}=I_{{K_n}\symbol{92} a}+ I_{{K_n}\symbol{92} \{b,d\}}$ are also splittings of $I_{K_n}$.}
\end{Remark}

\begin{Proposition} \label{directbasic}
    Let $n \geq 4$ be an integer and $I_{K_n}=I_{G_1}+I_{G_2}$ be a subgraph splitting of $I_{K_n}$. If $e$ is an edge of $K_n$ which does not belong to $G_1$ and $f$ is an edge of $K_n$ which does not belong to $G_2$, then the edges $e, f$ are adjacent in $K_n$. Moreover, if $h$ is another edge of $K_n$ which does not belong to $G_2$, then the edges $f, h$ are not adjacent in $K_n$.
\end{Proposition}
{\em \noindent Proof.} Suppose that the edges $e$, $f$ are not adjacent. Let $e=\{v_{i},v_{j}\}$, $f=\{v_{k},v_{l}\}$ and $e \cap f=\emptyset$. Let $\{g_{1},\ldots,g_{r}\}$ be a system of binomial generators of $I_{G_1}$ and $\{h_{1},\ldots,h_{s}\}$ be a system of binomial generators of $I_{G_2}$, then $\{g_{1},\ldots,g_{r},h_{1},\ldots,h_{s}\}$ is a generating set of $I_G$. Therefore we can find a minimal system of generators $V$ of $I_{K_n}$ formed by binomials belonging to either $I_{G_1}$ or $I_{G_2}$.
By Remark \ref{Indispe} the monomial $ef=\epsilon_{ij}\epsilon_{kl}$ is indispensable of $I_{K_n}$, so  it is a monomial term in a binomial $B_w$ of $I_{G_1}$ or $I_{G_2}$. Since $e$ is not an edge of $G_1$, we have that $B_{w} \notin I_{G_1}$. But $f$ is not an edge of $G_2$, so $B_{w} \notin I_{G_2}$. Thus $B_{w}$ does not belong to $I_{G_1}$ or $I_{G_2}$, a contradiction. Consequently, the edges $e$, $f$ are adjacent.\\
 Let $e=\epsilon_{ij}$ and $f=\epsilon_{jk}$. Suppose that $h$ is another edge of $K_n$ which does not belong to $G_2$. Since $e$ does not belong to $G_1$, we have, from the first part of the Proposition, that the edges $e$, $h$ are adjacent. We claim that $f$ and $h$ are not adjacent. Suppose that $f$, $h$ are adjacent. Since $h$ is adjacent to both $e$ and $f$, there are two cases for the edge $h$, namely either $h=\epsilon_{ik}$ or $h=\epsilon_{jl}$ for an index $l$ different than $i,j,k$. We will arrive at a contradiction when $h=\epsilon_{ik}$, and similarly, when $h=\epsilon_{jl}$, one also arrives at a contradiction. 
 
 Let $l$ be an index different than $i,j,k$. In a previous step we found that there is a minimal system of generators $V$ formed by binomials belonging to either $I_{G_1}$ or $I_{G_2}$. Moreover, $V$ must contain exactly two of the following three binomials $\epsilon_{ij}\epsilon_{kl}-\epsilon_{il}\epsilon_{jk}=e\epsilon_{kl}-\epsilon_{il}f, \epsilon_{ik}\epsilon_{jl}-\epsilon_{il}\epsilon_{jk}=h\epsilon_{jl}-\epsilon_{il}f,
\epsilon_{ij}\epsilon_{kl}-\epsilon_{ik}\epsilon_{jl}=e\epsilon_{kl}-h\epsilon_{jl}$.
But none of them belongs to $I_{G_2}$, since $f,h \not \in E(G_2)$, while $e\epsilon_{kl}-\epsilon_{il}f$ and $e\epsilon_{kl}-h\epsilon_{jl}$ do not belong to $I_{G_1}$, since $e \not \in E(G_1)$. A contradiction.

Consequently, the edges $f, h$ are not adjacent. \hfill $\square$

\begin{Remark} {\rm From Proposition \ref{directbasic} we deduce that any of the graphs $G_1, G_2$ contains all edges of $K_n$ except at most two. Suppose not, and let $e=\epsilon_{ij}\not \in E(G_1)$, $\{f_{1},f_{2},f_{3}\} \subset E(K_n)$ such that $f_l \not \in E(G_2)$, for every $1 \leq l \leq 3$. By Proposition \ref{directbasic}, any of the edges $f_1, f_2, f_3$ is adjacent to $e=\epsilon_{ij}$, so $v_i$ is a vertex of the edges $f_1, f_2, f_3$ or $v_j$ is a vertex of the edges $f_1, f_2, f_3$. Thus either $v_i$ or $v_j$ is a vertex of at least two of the edges $f_1, f_2, f_3$, and therefore two of the edges $f_1, f_2, f_3$ are adjacent, a contradiction to the second part of Proposition \ref{directbasic}.}
\end{Remark}

\begin{Proposition} \label{all splittings} Let $n \geq 4$ be an integer and $I_{K_n}=I_{G_1}+I_{G_2}$ be a subgraph splitting of $I_{K_n}$. Then there exists a cycle $w=(a,b,c,d)$ of length 4 in $K_n$, where $a, b, c, d \in E(G)$ such that 
\begin{enumerate}
    \item $G_1={K_n}\symbol{92} a$ and $G_2={K_n}\symbol{92} b$ or
    \item $G_1={K_n}\symbol{92} \{a,c\}$ and $G_2={K_n}\symbol{92} b$ or 
    \item $G_1={K_n}\symbol{92} a$ and $G_2={K_n}\symbol{92} \{b,d\}$ or 
    \item $G_1={K_n}\symbol{92} \{a,c\}$ and $G_2={K_n}\symbol{92} \{b,d\}$.
\end{enumerate}
\end{Proposition}
{\em \noindent Proof.} Since $I_{G_1} \subsetneqq I_{K_n}$, there exists an edge $a$ of $K_n$ such that $a$ is not an edge of $G_1$. Thus $I_{G_1} \subset I_{{K_n}\symbol{92} a}$. Since $I_{G_2} \subsetneqq I_{K_n}$, there exists an edge $b$ of $K_n$ such that $b$ is not an edge of $G_2$. Thus $I_{G_2} \subset I_{{K_n}\symbol{92} b}$. By Proposition \ref{directbasic}, the edges $a,b$ are adjacent. We distinguish the following cases:
\begin{enumerate} \item $G_1={K_n}\symbol{92} a$ and $G_2={K_n}\symbol{92} b$. Since $a,b$ are adjacent edges in the complete graph $K_n$ with $n \geq 4$ vertices, there exists a cycle $w$ of length 4 in $K_n$ with two consecutive edges $a, b$. 
\item $G_1 \subsetneqq {K_n}\symbol{92} a$ and $G_2 = {K_n}\symbol{92} b$ or $G_1={K_n}\symbol{92} a$ and $G_2={K_n}\symbol{92} \{b,d\}$. Suppose that there exists an edge $c$ of ${K_n}\symbol{92} a$ such that $c$ is not an edge of $G_1$. By Proposition \ref{directbasic}, the edges $a,c$ do not share a common vertex and they are adjacent to $b$. Furthermore, $G_1$ contains all edges of $K_n$ except at most two. Thus $G_1={K_n}\symbol{92} \{a,c\}$ and there exists a cycle $w$ of length 4 in $K_n$ with three consecutive edges $a, b, c$. The case (3), namely $G_1={K_n}\symbol{92} a$ and $G_2={K_n}\symbol{92} \{b,d\}$, is similar and therefore it is omitted.
\item[(3)] $G_1 \subsetneqq {K_n}\symbol{92} a$ and $G_2 \subsetneqq {K_n}\symbol{92} b$. Then there exists an edge $c$ of ${K_n}\symbol{92} a$ and an edge $d$ of ${K_n}\symbol{92} b$ such that $c$ is not an edge of $G_1$ and $d$ is not an edge of $G_2$. By Proposition \ref{directbasic}, the edges $a,c$ do not share a common vertex and they are adjacent to both $b, d$. Additionally $G_1$ contains all edges of $K_n$ except at most two. By the same proposition, the edges $b,d$ do not share a common vertex and they are adjacent to both $a, c$. Furthermore $G_2$ contains all edges of $K_n$ except at most two. Thus $w=(a,b,c,d)$ is a cycle in $K_n$, $G_1={K_n}\symbol{92} \{a,c\}$ and $G_2={K_n}\symbol{92} \{b,d\}$. \hfill $\square$
\end{enumerate}

 \begin{Theorem} \label{MinimalSplitting} Let $n \geq 4$ be an integer. Then $I_{K_n}$ has a minimal splitting if and only if $4\leq n \leq 5$.
 \end{Theorem}
 {\em \noindent Proof.} Suppose first that $n=4$. Then $I_{K_4}=I_{K_{4} \setminus \{\epsilon_{12}, \epsilon_{34}\}}+I_{K_{4} \setminus \{\epsilon_{14}, \epsilon_{23}\}}$ is a minimal splitting of $I_{K_4}$ by Example \ref{MinimalK_4}. Suppose now that $n=5$ and let $\{v_{1},\ldots,v_{5}\}$ be the vertex set of $K_5$. Let $G_1=K_{5} \symbol{92} \{\{v_{1},v_{2}\}, \{v_{3},v_{4}\}\}$ and $G_2=K_{5} \symbol{92} \{\{v_{1},v_{4}\}, \{v_{2},v_{3}\}\}$ be subgraphs of $K_5$. Then $S=\{\epsilon_{13}\epsilon_{24}-\epsilon_{14}\epsilon_{23}, \epsilon_{14}\epsilon_{25}-\epsilon_{15}\epsilon_{24}, \epsilon_{23}\epsilon_{45}-\epsilon_{24}\epsilon_{35},\epsilon_{13}\epsilon_{25}-\epsilon_{15}\epsilon_{23},\epsilon_{13}\epsilon_{45}-\epsilon_{14}\epsilon_{35}\}$ is a minimal generating set of $I_{G_1}$ and $T=\{\epsilon_{12}\epsilon_{34}-\epsilon_{13}\epsilon_{24}, \epsilon_{24}\epsilon_{35}-\epsilon_{25}\epsilon_{34}, \epsilon_{13}\epsilon_{45}-\epsilon_{15}\epsilon_{34},\epsilon_{12}\epsilon_{45}-\epsilon_{15}\epsilon_{24}, \epsilon_{12}\epsilon_{35}-\epsilon_{13}\epsilon_{25}\}$ is a minimal generating set of $I_{G_2}$. Also for any binomial $B \in S$ we have that neither $B$ nor $-B$ belongs to $T$, while for any binomial $B' \in T$ we have that neither $B'$ nor $-B'$ belongs to $S$. Moreover $S \cup T$ is a minimal generating set of $I_{K_5}$, and therefore $I_{K_5}=I_{G_1}+I_{G_2}$ is a minimal splitting of $I_{K_5}$. 

Finally assume that $n\ge 6$ and let $I_{K_n}=I_{G_1}+I_{G_2}$ be a minimal splitting of $I_{K_n}$. Then there exists a cycle $w=(a,b,c,d)$ of $K_n$ such that $G_1$ and $G_2$ are of one of the four types of Proposition \ref{all splittings}. Without loss of generality we can assume that $a=\{v_{1},v_{2}\}$, $b=\{v_{2},v_{3}\}$, $c=\{v_{3},v_{4}\}$ and $d=\{v_{1},v_{4}\}$. Since $n\geq 6$, the graph $K_n$ has at least two more vertices, say $v_5$ and $v_6$. Notice that the complete subgraph of $K_n$ on the vertex set $\{v_1, v_3, v_5, v_6\}$ is also a subgraph of both $G_1$ and $G_2$. Let $S$ be a minimal system of binomial generators of $I_{G_1}$ and $T$ be a minimal system of binomial generators of $I_{G_2}$. Both $S$ and $T$ must contain exactly two of the binomials $\epsilon_{13}\epsilon_{56}-\epsilon_{15}\epsilon_{36}, \epsilon_{13}\epsilon_{56}-\epsilon_{16}\epsilon_{35}, \epsilon_{15}\epsilon_{36}-\epsilon_{16}\epsilon_{35}$. Thus $S$ and $T$ have at least one minimal generator in common which contradicts the fact that $I_{K_n}=I_{G_1}+I_{G_2}$ is a minimal splitting of $I_{K_n}$. Consequently, for $n \geq 6$ the ideal $I_{K_n}$ has no minimal splitting. \hfill $\square$\\

\section{Reduced splittings}

In this section, we introduce reduced splittings of toric ideals of graphs and show that every minimal splitting of the toric ideal of a graph is also a reduced splitting.

\begin{Definition} We say that the subgraph splitting $I_G=I_{G_1}+I_{G_2}$ of $I_G$ is reduced if for any subgraphs $H_1$ of $G_1$ and $H_2$ of $G_2$ with $I_G=I_{H_1}+I_{H_2}$ it holds that $I_{H_1}= I_{G_1}$ and $I_{H_2}= I_{G_2}.$

 \end{Definition}
A basic step towards determining all subgraph splittings of the toric ideal of a graph is to find its reduced splittings. All other subgraph splittings are found from a reduced splitting $I_G=I_{G_1}+I_{G_2}$ by adding edges to one of $G_1$, $G_2$ or both to get graphs $G_1', G_2'$, as long as $I_G=I_{G_1'}+I_{G_2'}$ is a splitting.

\begin{Remark} {\rm The reduced splittings of $I_{K_n}$ are those of the last type in Proposition \ref{all splittings}, namely $I_{K_n}=I_{G_1}+ I_{G_2}$ where $G_1={K_n}\symbol{92} \{a,c\}$ and $G_2={K_n}\symbol{92} \{b,d\}$. The other two types in Proposition \ref{all splittings} can be taken from the reduced splittings, by adding edges.}
\end{Remark}

\begin{Proposition} If $I_G$ is subgraph splittable, then it has at least one reduced splitting.
\end{Proposition}
{\em \noindent Proof.} Any subgraph splitting $I_G=I_{G_1'}+I_{G_2'}$ of $I_G$ is either reduced or there exist a splitting $I_G=I_{G_1}+I_{G_2}$ such that $G_{1}$ is a proper subgraph of $G_1'$ or/and $G_{2}$ is a proper subgraph of $G_2'$. In the latter case, $G_{1}$ has fewer edges than $G_1'$ or/and $G_{2}$ has fewer edges than $G_2'$. This procedure cannot be repeated indefinitely, since the number of edges of $G$ is finite.  \hfill $\square$\\

 To understand the structure of reduced splittings one has to generalize first the notion of edge splitting by replacing the edge with a set of edges.

Let $S=\{B_{w_1}, B_{w_2}, \ldots, B_{w_r}\}$ be a minimal system of binomial generators of $I_G$. Given a set $F\subset E(G)$, we define $G_S^F=\bigcup_{e\in F}G_S^e$, so $I_{G_S^F}=\sum_{e \in F}I_{G_S^e}$ and also $I_G=I_{G_S^F}+I_{G\symbol{92}F}$, by using similar arguments as in the proof of Theorem \ref{sum1}. 

Of particular interest is the case that $F$ is the set of all edges of $G$ having a common vertex $v$. Given a vertex $v$ of $G$, we let $G-v$ be the subgraph of $G$ obtained by deleting the vertex $v$. We denote by $G_{S}^{v}$ the subgraph of $G$ with edges $$E(G_{S}^{v})=\bigcup_{1 \leq i \leq r \ \textrm{and} \ v \in V(w_i)} E(w_i).$$ It holds that $I_G=I_{G_S^{v}}+I_{G-v}$.

The next theorem asserts that the reduced splittings of $I_G$ are always in the form $I_G=I_{G_S^F}+I_{G\symbol{92}F}$, for suitable sets $S$ and $F$.
\begin{Theorem}
    Let $I_G=I_{G_1}+I_{G_2}$ be a reduced splitting of $I_G$. Then there exist a set $F\subset E(G)$ and a minimal system of binomial generators $S$ of $I_G$ such that $I_{G_1}=I_{G_S^F}$ and $I_{G_2}=I_{G\symbol{92}F}$.
\end{Theorem}
{\em \noindent Proof.} Let $I_G=I_{G_1}+I_{G_2}$ be a reduced splitting of $I_G$ and set $F=G \symbol{92} G_2$. Let $\{f_{1},\ldots,f_{s}\}$ be a system of binomial generators of $I_{G_1}$, $\{g_{1},\ldots,g_{t}\}$ be a system of binomial generators of $I_{G_2}$ and $S=\{B_{w_1}, \ldots,B_{w_r}\} \subset \{f_{1},\ldots,f_{s},g_{1},\ldots,g_{t}\}$ be a minimal system of binomial generators of $I_G$ as in the proof of Theorem \ref{basicsplit}. Then $G\symbol{92}F=G_{2}$, so $I_{G\symbol{92}F}=I_{G_2}$, and $I_{G_{S}^{F}} \subset I_{G_1}$ since $I_{G_{S}^e} \subset I_{G_1}$ for each $e\in F$ from the proof of Theorem \ref{basicsplit}. But $I_G=I_{G_1}+I_{G_2}$ is a reduced splitting of $I_G$ and also $I_G=I_{G_S^F}+I_{G_2}$, since $I_G=I_{G_S^F}+I_{G\symbol{92}F}$ and $I_{G\symbol{92}F}=I_{G_2}$, with $I_{G_S^F}\subset I_{G_1}$, therefore $I_{G_S^F}= I_{G_1}$.
\hfill $\square$

\begin{Remark} {\em A reduced splitting $I_G=I_{G_1}+I_{G_2}$ can be also written in the form $I_G=I_{G\symbol{92}F}+I_{G_S^F}$ where $F=G \symbol{92} G_1$, $I_{G_1}=I_{G\symbol{92}F}$ and $I_{G_2}=I_{G_S^F}$.}
\end{Remark}

\begin{Example} {\em Let $G$ be the bipartite graph consisting of four 4-cycles $w_1, w_2, w_3, w_4$ in a row, i.e. for $i<j$ it holds that $E(w_i)\cap E(w_j)=\emptyset$ except if $j=i+1$
in which case they have one edge in common. 
\begin{figure}[ht]
\begin{center}
\includegraphics[scale=0.35]{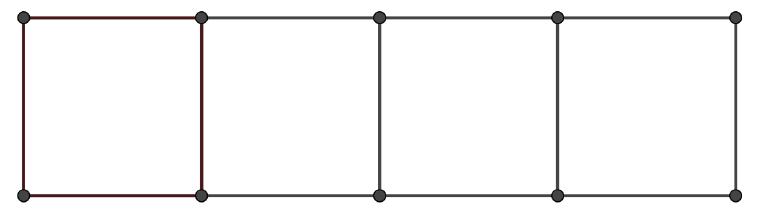}
\caption{The graph $G$ with four 4-cycles in a row}
\label{graph4}
\end{center}
\end{figure}

The ideal $I_{G}$ has a unique minimal system of binomial generators consisting of the binomials $B_{w_1}, B_{w_2}, B_{w_3}, B_{w_4}$. Then there are 25 different splittings. More precisely, four minimal and reduced splittings in the form $I_G= \langle B_{w_i} \rangle+ \langle B_{w_j}, B_{w_k}, B_{w_l} \rangle$, where $\{i,j,k,l\}=\{1,2,3,4\}$. Also three minimal and reduced splittings in the form $I_G=\langle B_{w_i}, B_{w_j} \rangle+ \langle B_{w_k}, B_{w_l} \rangle$. Finally eighteen  non-minimal and non-reduced splittings, twelve in the form $I_G=\langle B_{w_i}, B_{w_j} \rangle+ \langle B_{w_j}, B_{w_k}, B_{w_l} \rangle$ and six in the form $I_G=\langle B_{w_i}, B_{w_k},  B_{w_l} \rangle+ \langle B_{w_j}, B_{w_k}, B_{w_l} \rangle$.}
\end{Example}

The next theorem asserts that minimal splittings are always reduced.

\begin{Theorem} \label{ReducedMinimal} Every minimal splitting of $I_G$ is also a reduced splitting. 
    
\end{Theorem}
{\em \noindent Proof.} Let $I_G=I_{G_1}+I_{G_2}$ be a minimal splitting which is not reduced. Then there exist  $I_{G_1'}\subset I_{G_1}$ and $I_{G_2'}\subset I_{G_2}$ such that $I_G=I_{G_1'}+I_{G_2'}$, where $I_{G_{1}'}$ is a proper subset of $I_{G_1}$ or/and $I_{G_{2}'}$ is a proper subset of $I_{G_2}$. Suppose that for instance $I_{G_{1}'}$ is a proper subset of $I_{G_1}$. Let $\{B_{w_1},\dots, B_{w_s}\}$ and $\{B_{w_{s+1}}, \dots, B_{w_l}\}$ be minimal systems of binomial generators of the ideals $I_{G_1}$ and $I_{G_2}$, respectively. Since $I_G=I_{G_1}+I_{G_2}$ is a minimal splitting, the ideal $I_G$ is minimally generated by the set $\{B_{w_1},\dots, B_{w_s}, B_{w_{s+1}}, \dots, B_{w_l}\}$. Let $\{f_1, \dots, f_t\}$ be a system of binomial generators of $I_{G_1'}$. Then from the equality $I_G=I_{G_1'}+I_{G_2'}$ we have that $I_G=I_{G_1'}+I_{G_2}$, since $G_2'\subset G_2$. Thus there exists a set $\{B_{w_1'},\dots, B_{w_s'}\} \subset \{f_1, \dots, f_t\}$ such that $\{B_{w_1'},\dots, B_{w_s'}, B_{w_{s+1}}, \dots, B_{w_l}\}$ is a minimal system of generators of $I_G$, since toric ideals of graphs are homogeneous and therefore any two minimal systems of generators have the same cardinality. Then, after reordering $B_{w_1'},\dots, B_{w_s'}$ if necessary, we can assume that $B_{w_j'}=B_{w_j}$ if $B_{w_j}$ is indispensable, or that $B_{w_j'}$ and $B_{w_j}$ are $F_4$-equivalent if $B_{w_j}$ is dispensable. Recall that two primitive walks $\gamma$, $\gamma'$ are $F_4$-equivalent if either $\gamma=\gamma'$ or there exists a series of walks $\gamma_{1}=\gamma, \gamma_{2},\ldots,\gamma_{r-1}, \gamma_{r}=\gamma'$ such that $\gamma_{i}$ and $\gamma_{i+1}$ differ by an $F_4$, where $1 \leq i \leq r-1$, see \cite[Section 4]{RTT}. Since $G_{1}'$ is a proper subgraph of $G_1$, there exists an edge $e\in E(w_i) \subset E(G_1)$ which is not in $E(G_1')$, but $e\not \in E(w_i')$ for at least one index $1\leq i\leq s$. Then $B_{w_i}$ is dispensable, thus $w_i, w_i'$ are $F_4$-equivalent and $e$ belongs to a common $F_4$ of both $w_i, w_i'$. Suppose that the edges of the $F_4$ belonging to $w_i$ are $e,f$ and to $w_i'$ are $a,b$, thus $F_{4}=(e,a,f,b)$. Consider the binomial $ef-ab \in I_G$, we have that all $e,f,a,b$ are edges of $G_1$, since $G_1' \subset G_1$, and therefore $ef-ab\in I_{G_1}$.

We distinguish the following cases. \begin{enumerate}
\item $ef-ab$ is an indispensable binomial of $I_G$. The set $\{B_{w_1'},\dots, B_{w_s'}, B_{w_{s+1}}, \dots, B_{w_l}\}$ is a minimal system of generators of $I_G$ and $ef-ab$ is indispensable of $I_G$. So $ef-ab$ is one of the binomials $B_{w_i'}, 1 \leq k \leq s$ or $B_{w_k}, s+1\leq k\leq l$. But $ef-ab\not \in I_{G_1'}$ since $e$ is not an edge of $G_1'$, thus $ef-ab=B_{w_k}$ for an index $s+1\leq k\leq l$. Since the binomial $ef-ab$ is indispensable of $I_G$, the $A_G$-fiber of ${\rm deg}_{A_{G}}(ef)$ has only two elements, namely $ef$ and $ab$, see \cite{ChKTh}. Then the $A_{G_1}$-fiber of ${\rm deg}_{A_{G_1}}(ef)$ has at most two elements; in fact, it has exactly two with no common factor other than 1 since $ef-ab\in I_{G_1}$. Thus $ef-ab$ is indispensable of $I_{G_1}$, and therefore $ef-ab=B_{w_q}$ for an index $1\leq q \leq s$ a contradiction to the hypothesis that $I_G=I_{G_1}+I_{G_2}$ is a minimal splitting.
\item $ef-ab$ is not an indispensable binomial of $I_G$. Then there is a binomial $ef-cd$ in $I_G$. In this case, the $A_G$-fiber of ${\rm deg}_{A_G}(ef)$ corresponds to a subgraph of $G$ homomorphic to $K_4$ and it consists of exactly three monomials, namely $ef$, $ab$, and $cd$. Every minimal system of generators of $I_G$ should contain exactly two of the binomials $ef-ab, ef-cd, ab-cd$. There are two subcases.
\begin{enumerate}
 \item[(i)] Both edges $c, d$ belong to $G_1$. Then all binomials $ef-ab, ef-cd, ab-cd$ belong to $I_{G_1}$, so any minimal system of binomial generators of $I_{G_1}$ should contain exactly two of them. Since $I_G=I_{G_1}+I_{G_2}$ is a minimal splitting, the ideal $I_{G_2}$ cannot contain any of the above three binomials. But $I_G=I_{G_1'}+I_{G_2}$ is also a splitting of $I_G$ and $e$ is not an edge of $G_1'$, therefore only $ab-cd$ can be an element of $I_{G_1'}$. Then $\{B_{w_1'},\dots, B_{w_s'}, B_{w_{s+1}}, \dots, B_{w_l}\}$ is a minimal generating set of $I_G$, which contains at most one of the binomials $ef-ab, ef-cd, ab-cd$, a contradiction.
\item[(ii)] At least one of $c, d$ does not belong to $G_1$. Then the binomial $ef-ab$ is indispensable of $I_{G_{1}}$, since the $A_{G_1}$-fiber of ${\rm deg}_{A_{G_1}}(ef)$ has only two elements, namely $ef$ and $ab$, with no common factor other than 1. Thus the set $\{B_{w_1},\dots, B_{w_s}\}$ contains the binomial $ef-ab$. The set $\{B_{w_1},\dots, B_{w_s}, B_{w_{s+1}}, \dots, B_{w_l}\}$ is a minimal system of generators of $I_G$, so $\{B_{w_{s+1}}, \dots, B_{w_l}\}$ contains exactly one of the binomials $ef-cd$ or $ab-cd$ and does not contain $ef-ab$.
But $I_G=I_{G_1'}+I_{G_2}$ is a splitting of $I_G$ and none of the binomials $ef-ab, ef-cd, ab-cd$ belongs to $I_{G_1'}$, since $e$ and at least one of $c, d$ does not belong to $E(G_1')$. 
Thus the set $\{B_{w_{s+1}}, \dots, B_{w_l}\}$ contains exactly two of the binomials $ef-ab, ef-cd, ab-cd$, a contradiction.
\end{enumerate}
\end{enumerate}
In all cases we reach a contradiction, so every minimal splitting is reduced. 
\hfill $\square$

\begin{Remark} {\rm The converse of Theorem \ref{ReducedMinimal} is not true. By Theorem \ref{MinimalSplitting}, $I_{K_n}$ does not have a minimal splitting for $n \geq 6$. But from Proposition \ref{all splittings} any splitting $I_G=I_{{K_n}\symbol{92} \{a,c\}}+I_{{K_n}\symbol{92} \{b,d\}}$ is reduced.}
\end{Remark}

\begin{Remark} \label{Square} {\rm 
In \cite[Question 5.1]{FHKT} G. Favacchio, J. Hofscheier, G. Keiper, and A. Van Tuyl pose the following question: For what graphs $G$ can we find $G_1$ and $G_2$ so that their respective toric ideals satisfy $I_G=I_{G_1}+I_{G_2}$? We answered this question in the case that $G_1, G_2$ are subgraphs of $G$. The general question seems difficult, since one has information from $G$ only about the edges and possible closed walks of $G_1, G_2$, but not about their set of vertices. }
\end{Remark}

The next example shows that there may be splittings $I_G=I_{G_1}+I_{G_2}$ such that one of $G_1, G_2$ or both are not isomorphic to any subgraph of $G$. 
\begin{Example} \label{Nosubgraphsplittins} {\rm

Let $G$ be the grid graph of Figure \ref{square} with 9 vertices and 12 edges.
\begin{figure}[ht]
\begin{center}
\includegraphics[scale=0.32]{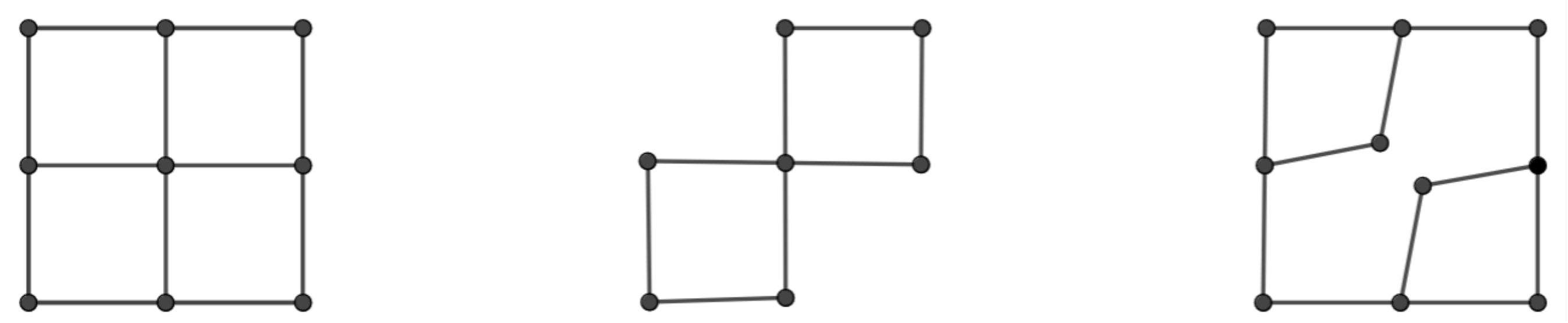}
\caption{The graphs $G, G_1$ and $G_2$ of a non subgraph splitting $I_{G}=I_{G_1}+I_{G_2}$.}
\label{square}
\end{center}
\end{figure}
The graph $G$ has four cycles of length 4 with no chord, namely $w_{1}$, $w_{2}$, $w_{3}$ and $w_{4}$, and a cycle of length 8 with no chord, namely $w_{5}$. Then the ideal $I_G$ is generated by the binomials $B_{w_i}, 1 \leq i \leq 5$. We consider the subgraph $G_1$ of $G$ consisting of the vertices and edges of the walks $w_1, w_3$. Let $G_2$ be the third graph of Figure \ref{square} with 10 vertices and 12 edges. We may think $G$ as the graph obtained from $G_2$ by identifying the two vertices in the center of the figure of $G_2$.
In this way there is a correspondence between the edges of $G$ and $G_2$ and thus the  $I_G, I_{G_1}, I_{G_2}$ are ideals in the same polynomial ring $K[e_1, e_2, \dots , e_{12}].$
The ideal $I_{G_1}$ is generated by $B_{w_1}, B_{w_3}$ and the ideal $I_{G_2}$ is generated by $B_{w_2}, B_{w_4}, B_{w_5}, B_{w_6},  B_{w_7}, B_{w_8}$, where $w_i,\ 5 \leq i \leq 8$, are closed walks of length 8 with no chord. Then $I_{G}=I_{G_1}+I_{G_2}$ is a splitting of $I_{G}$ and $G_2$ is not isomorphic to any subgraph of $G$, since it has 10 vertices. Thus there exist splittings of graphs that are not subgraph splittings.}
\end{Example}
\noindent{\bf Acknowledgements} The authors would like to thank the two anonymous reviewers for their positive and constructive feedback
and for providing valuable suggestions that have enhanced the clarity, readability and presentation of the paper.\\

\noindent \textbf{Funding} The authors did not receive support from any organization for the submitted work.\\

\noindent \textbf{Data availability} Data sharing not applicable to this article as no datasets were generated or analysed during the current study.\\

\noindent \textbf{Declaration}\\

\noindent \textbf{Conflict of interest}  The authors declare that there are no conflict of interest that are directly or indirectly related to the work submitted for publication.

\end{document}